\documentclass[a4paper,10pt]{article} 
 
\usepackage{fancybox} 
\usepackage{textcomp}
\usepackage{pifont}
\usepackage{siunitx}
\usepackage[stable]{footmisc}
\usepackage{graphicx}
\usepackage{caption}
\graphicspath{{pictures/}}
\DeclareGraphicsExtensions{.pdf,.png,.jpg}
\usepackage{comment}
\usepackage[russian, english]{babel}
\usepackage{fancybox} 
\usepackage{amssymb} 
\usepackage{amsmath} 
\usepackage[applemac]{inputenc}
\usepackage{amsfonts}
\usepackage{amsthm}
\usepackage{color}
\usepackage{mathrsfs} 
\usepackage{dsfont}
\usepackage{titling}
\usepackage{hyperref}
\hypersetup{
	colorlinks=true,
	linkcolor=blue,
	filecolor=magenta,      
	urlcolor=cyan,
}
\usepackage{mathtools}                                                      
\mathtoolsset{showonlyrefs=true}      

\usepackage{color}

\newcommand\blfootnote[1]{%
	\begingroup
	\renewcommand\thefootnote{}\footnote{#1}%
	\addtocounter{footnote}{-1}%
	\endgroup
}

\newtheorem{theorem}{Theorem}[section]
\newtheorem{proposition}[theorem]{Proposition}
\newtheorem{corollary}[theorem]{Corollary}
\newtheorem{lemma}[theorem]{Lemma}

\newtheorem*{remark}{Remark}
\newtheorem*{remarks}{Remarks}

\def \R {{\mathbb {R}}}
\def \N {{\mathbb {N}}}

\def \phi {{\varphi}}

\def \tilde {\widetilde}

\thanksmarkseries{arabic}
\numberwithin{equation}{section}

\title{\textbf{$L_1$ and $L_{\infty}$ stability of transition densities of perturbed diffusions}\protect\blfootnote{\small The study has been funded by the Russian Science Foundation (project \textnumero \ 20-11-20119)}}
\author{\textbf{I.Bitter}\thanks{Laboratory of Stochastic Analysis, HSE University,
		Pokrovsky Blvd, 11, Moscow, Russian Federation. ilya.bitter@yandex.ru}, \textbf{V. Konakov}\thanks{Laboratory of Stochastic Analysis, HSE University,
		Pokrovsky Blvd, 11, Moscow, Russian Federation. vkonakov@hse.ru}}

\begin{document}
\maketitle 

{\small{\textbf{MSC \textcolor{black}{2020}:} Primary: 60H10, 60H30; Secondary: 35K10.}}

\begin{abstract}
	In this paper, we derive a stability result for $L_1$ and $L_{\infty}$ perturbations of diffusions under weak regularity conditions on the coefficients. In particular, in contrast to \cite{KKM, KM}, the drift terms we consider can be unbounded with at most linear
	growth, and we do not require uniform convergence of perturbed diffusions.
	Instead, we require a weaker convergence condition in a special metric introduced
	in this paper, related to the Holder norm of the diffusion matrix differences. Our
	approach is based on a special version of the McKean-Singer parametrix
	expansion.
\end{abstract}

\tableofcontents	
\section{Introduction} \label{Intro}

\subsection{Setting}
	For a fixed deterministic horizon $T > 0,$ let us consider the following 
$d$ - dimensional, non-degenerate Ito's diffusion
\begin{equation}\label{SDE_0}
d X_{t}=b\left(t, X_{t}\right) d t+\sigma\left(t, X_{t}\right) d W_{t}, \quad t \in[0, T].
\end{equation}
$W_t$ stands for a  $d$ - dimensional Brownian motion on some filtered complete probability space $(\Omega,\mathcal F,(\mathcal F_t)_{t \geq 0},\mathbb P) $ satisfying the usual conditions, and the coefficients are assumed to be measurable, rough in time, and H\"older in space, with potentially unbonded drift. Also $a(t, x):=\sigma \sigma^{*}(t, x)$ is assumed to be uniformly elliptic. The infinitesimal generator of \eqref{SDE_0} at time $u$ and for all $\varphi \in C_{0}^{2}\left(\mathbb{R}^{d}, \mathbb{R}\right), z \in \mathbb{R}^{d}$, is defined as follows:
\begin{align}\label{generator}
L_{u} \varphi(z)=\frac{1}{2} \operatorname{Tr}\left(\sigma \sigma^{*}(u, z) D_{z}^{2} \varphi(z)\right)+\left\langle b(u, z), D_{z} \varphi(z)\right\rangle =& \\ = \frac{1}{2} \sum_{i, j=1}^{d}a_{i j}(u, z) \frac{\partial^{2} \phi(u, z)}{\partial z_{i} \partial z_{j}}+\sum_{i=1}^{d}b_{i}(u, z) \frac{\partial \phi(u, z)}{\partial z_{i}}.
\end{align}

We now introduce, for a given parameter $\varepsilon > 0,$ a perturbed version of \eqref{SDE_0} with dynamics
\begin{equation}\label{SDE_1}
d X_{t}^{\varepsilon}=b_{\varepsilon}\left(t, X_{t}^{\varepsilon}\right) d t+\sigma_{\varepsilon}\left(t, X_{t}^{\varepsilon}\right) d W_{t}, \quad t \in[0, T],
\end{equation} where the coefficients $b_{\varepsilon}, \sigma_{\varepsilon}$ satisfy at least the same assumptions as $b, \sigma$ and are in some sense meant to be close to $b, \sigma$ when $\varepsilon$ is small. When both coefficients $b, \sigma$ are bounded and Holder continuous, it is well known that there exists a unique weak solution to \eqref{SDE_0}, which has a density \cite{fried, kalash, KKM}. Moreover, combining the parametrix method to get the upper bound \cite{KKM, deck}, and the chaining method to obtain the lower bound \cite{bass}, it can be shown \cite{Ar1, Ar2} that the transition density $p(0,x,t,y)$ satisfies the two sided Gaussian bounds  
$$C^{-1} g_{\lambda^{-1}}(t, x-y) \leq p(0, x, t, y) \leq C g_{\lambda}(t, x-y),$$ where $$g_{\lambda}(t, x)=t^{-d / 2} \exp \left(-\frac{\lambda|x|^{2}}{t}\right), \lambda \in (0,1], \quad t>0.$$ 

Such methods were successfully developed for more general cases, namely, for operators satisfying strong H\"ormander conditions and Kolmogorov operators with linear drifts \cite{KMM, pol, DFP}.  When the drift is unbounded and non-linear, fewer results are available. To obtain an upper bound, we need to control the terms of the parametrix series, and in the case of unbounded drift it becomes a delicate problem. For the drifts with sublinear growth, namely, 
$|b(t, x)| \leq C\left(|x|^{\beta}+1\right), \beta \in(0,1),$ the generalization of the parametrix method was obtained in \cite{deck}, but the method developed there fails for the drifts with linear growth.  It seems quite plausible that a linearly growing drift is exactly the boundary case, starting from which, it is necessary to introduce a forward flow corresponding to the transport of the initial condition or, equivalently, a backward flow corresponding to the transport of the terminal condition \cite{MPZ, DM, KMM}. For unbounded at most linearly growing drift, the upper bound for the transition density may be obtained by using the truncation method introduced in \cite{DM}, see also \cite{MPZ}. This method consists in viewing a Fleming logarithmic transformation (see \cite{DM, sheu}) of a transition density as a value function of a certain stochastic control problem \cite{cor}, the desired density estimates are obtained by choosing appropriate controls. In the case of unbounded drift, the truncation method allows one to obtain an upper bound for a transition density, but, in order to study the sensitivity of transition densities w.r.t. the perturbations of coefficients one needs to work with a complete parametrix series without truncation.                

Note also that the upper bounds for the total variation, entropy and Kantorovich distances between two solutions to Fokker-Planck-Kolmogorov equations with different diffusion matrices and drifts on $\mathbb{R}^{d} \times[0, T]$ with fixed $T > 0$ were obtained in \cite{bog-shap}, and upper bounds for the distances between stationary solutions to Fokker-Planck-Kolmogorov equations with different diffusion matrices and drifts were obtained in \cite{bog-shap-2}.

The aim of the present paper is to extend the stability results \cite{KKM} to the case of nonlinear unbounded drift with at most linear growth. Moreover, we replace the condition of $L_{\infty}$ - closeness by the weaker one, namely, by $L_1$ - closeness.  To this end, the complete parametrix series will be constructed, which, in contrast to \cite{deck}, will be convergent also for drifts with linear growth.                                        

Diffusions with the dynamics \eqref{SDE_0} and unbounded drifts appear in many applied problems. Important applications can, for instance, be found in mathematical finance. It is often very useful to know how a perturbation of the volatility $\sigma$ impacts the density, and therefore the associated option prices  \cite{CFP, BGM}. In the framework of parameter estimation, it can be useful, having at the hand estimators 
$\left(b_{\varepsilon}, \sigma_{\varepsilon}\right)$ of the true parameters $\left(b, \sigma\right)$ and some controls for the differences $\left|b - b_{\varepsilon}\right|, \left|\sigma - \sigma_{\varepsilon}\right|$ in a suitable sense, to quantify the difference $p - p_{\varepsilon}$ of the densities corresponding, respectively, to the dynamics with the estimated parameters and the one of the model. Another important application includes the case of mollification by spatial convolution. This specific kind of perturbation is useful for investigating the error between the densities of a non-degenerate diffusion of type \eqref{SDE_0} with H\"older coefficients (or with piecewise smooth bounded drift) and its Euler scheme. In this framework, some explicit convergence results can be found in \cite{KM} and for measurable bounded drift in \cite{BJ}. In summary, stability results can be useful in every applied field where the diffusion coefficients might be misspecified.

\bigskip

This paper is organized as follows. In the next subsection, we introduce our assumptions and state the main result. Section 2 contains an introduction to some crucial facts about parametrix expansion and deterministic flow associated with a drift of corresponding diffusion. In Section 3, we derive our main results in $L_1$ and $L_{\infty}$ cases, respectively. Additional proofs of technical lemmas are given in Appendix. In what follows, the letter $C$ stands for a positive constant, only depending on the quantities in the assumptions; in the context of the proofs, it's value may vary from line to line. 

\subsection{Assumptions and Main Result}\label{results}

In this paper, we study the difference between transition densities of the process \eqref{SDE_0} and its perturbed version \eqref{SDE_1}. A similar problem for processes with bounded drift was discussed in \cite{KKM}. In contrast, we consider a more general class of processes with growing drift. In our setting, we prove that the difference between $p$ and $p_{\varepsilon}$ admits an upper bound in terms of $L_1 - L_1$ norm: for probability measure $\mu$,  we denote  $$||f(x, y)||_{L_1^{\mu}\left(\mathbb{R}^d, L_1(\mathbb{R}^d)\right)} = \int_{\mathbb{R}^d}\int_{\mathbb{R}^d}|f(x, y)| dy \ \mu(dx).$$

For any Borel sets $A \in \mathcal{B}([t,s]), \ B \in \mathcal{B}(\mathbb{R}^d)$ and $x \in \mathbb{R}^d$, define the probability measures $$
\beta(A)=\int_{A} \mathfrak{B}\left(u ; 1, \frac{\gamma}{2}\right) d u, \ \bar{\lambda}(B)=\int_{B} \bar{p}(t, s, x, y) d y,
$$

and $L_{1}$-space of functions defined on $\left([t, s] \times \mathbb{R}^{d} \times \mathbb{R}^d, \beta \times \bar{\lambda} \times \mu \right)$
\begin{equation}
	||f(u,x,y)||_{L_{1}} =  
	\left[\int_{\mathbb{R}^d}\int_{R^{d}}\int_{t}^{s}\mathfrak{B}(u ; 1, \dfrac{\gamma}{2}) \left|f(u,x, y) \right|\bar{p}(t,s,x,y) \ du \  dy \ \mu(dx)\right].
\end{equation}
Here $\mathfrak{B}(u ; 1, \dfrac{\gamma}{2})$ is a beta-density on the interval $[t,s]$, $0 \leq t < s \leq T$, and $\bar{p}(t,s,x,y)$ is the transition density of some an auxiliary diffusion process that will be introduced in \eqref{Maj_dens}. Recall that  $$\mathfrak{B}(u; \alpha, \beta) = \dfrac{1}{B(\alpha, \beta)} (u - t)^{\alpha - 1} (s - u)^{\beta - 1} (s - t)^{1 - \alpha - \beta},\ u \in [t, s],$$ where $B(\alpha, \beta)$ is the Beta-function. It is natural in this case to consider the perturbation in $L_{1}$ sense and to use the beta-density as a weight function in the time variable. This choice is well-adapted to the beta function appearing in the upper bounds for the parametrix series terms.
\bigskip

Let us introduce the following assumptions. Below, the
parameter  $\varepsilon > 0$ is fixed and the constants appearing in the
assumptions do not depend on it. In what follows, it is important that we consider the problem globally, namely on the product space $[0, T] \times \mathbb{R}^d$.

\begin{itemize}
	\item[\textbf{(A1)}] (\textbf{Uniform Ellipticity.})
	Matrices $a = \sigma \sigma^{*}$ and $a_\varepsilon = \sigma_\varepsilon \sigma^{*}_\varepsilon$ are uniformly elliptic, i.e., there exists a positive constant $\Lambda \geq 1,$ such that $\forall (t, x) \in [0, T] \times \mathbb{R}^d$,
	\begin{align}
	&\Lambda^{-1}|w|^{2} \leq \langle a(t, x) w, w\rangle \leq \Lambda|w|^{2}, \quad \forall w \in \mathbb{R}^{d}. \\
	&\Lambda^{-1}|w|^{2} \leq \langle a_\varepsilon(t, x) w, w\rangle \leq \Lambda|w|^{2}, \quad \forall w \in \mathbb{R}^{d}.
	\label{assumption1}
	\end{align}
	
	This condition also implies the uniform ellipticity of matrices $a^{-1}$ and $a^{-1}_\varepsilon$.
	\item[\textbf{(A2)}] (\textbf{Regularity.})
	The diffusion coefficients and the drifts of \eqref{SDE_0} and \eqref{SDE_1} are H\"older and Lipschitz continuous  in space, respectively:
	for $\gamma \in(0,1], K > 0 \ \text{and} \ \forall(t, x),(s, y) \in [0, T] \times \mathbb{R}^d$, 
	\begin{equation}
	\left|\sigma(t, x) - \sigma(t, y)\right| + \left|\sigma_\varepsilon(t, x) - \sigma_\varepsilon(t, y)\right| \leq K\left|x-y \right|^{\gamma},
	\end{equation}
	\begin{equation}
	\left|b(t, x) - b(t, y)\right| + \left|b_\varepsilon(t, x) - b_\varepsilon(t, y)\right| \leq K\left|x-y\right|,
	\end{equation}
	\begin{equation}
	\left|b(t, \textbf{0})\right| + \left|b_\varepsilon(t, \textbf{0})\right| \leq K, \ \textbf{0} \in \mathbb{R}^d.
	\label{assumption2}
	\end{equation} 
Clearly, if $\sigma(t, \cdot)$ is bounded and H\"older continuous so is $a(t, \cdot)$. Also, the regularity conditions imply at most linear growth of the drift coefficients:
	\begin{align}
	&|b(t, x)| + |b_\varepsilon(t, x)| \leq K\left(1 + |x|\right) \quad \forall(t, x) \in [0, T] \times \mathbb{R}^d. \\
	\label{assumption3}
	\end{align}

	Thanks to the Lipschitz continuity and at most linear growth of the drift term in \eqref{SDE_0}, the following Cauchy problem has a unique solution for fixed $(s, y) \in [0, T] \times \mathbb{R}^d$, defined globally on $[0,s]$:

	\begin{equation}
	\dot{\theta}_{w, s}(y) = b(w, \theta_{w, s}(y)), \quad \theta_{s, s}(y) = y.
	\label{flow}
	\end{equation}
	
	Throughout the paper, we will call the deterministic flow solving ODE \eqref{flow} \textit{the flow associated with drift} $b$. In the same way, $\theta^{\varepsilon}_{w, s}(y)$ is the flow associated with $b_{\varepsilon}$.
	
	 \item[\textbf{(A3)}] (\textbf{Almost equivalence of the flows.})
	 
	 There exists a positive constant $C_{\theta}$ such that $\forall 0 \leq t < s \leq T$, uniformly in space,
	 
	 \begin{equation}
	 \int_t^s \left|b_{\varepsilon}(u, x) - b(u, x)\right| du \leq C_{\theta} \sqrt{s - t}.
	 \label{drift_diff}
	 \end{equation}
	 As a consequence, the flows $\theta_{w, s}(y)$ and $\theta^{\varepsilon}_{w, s}(y)$ associated with the drifts of \eqref{SDE_0} and \eqref{SDE_1} are almost equivalent:
	\begin{gather}
	\left|\theta^{\varepsilon}_{t, s}(y) - \theta_{t, s}(y)\right| = \left| \int_t^s b_{\varepsilon}\left(u, \theta^{\varepsilon}_{u, s}(y)\right) - b\left(u, \theta_{u, s}(y)\right) du \right|  \leq \\ \leq K\int_t^s \left|\theta^{\varepsilon}_{u, s}(y) - \theta_{u, s}(y) \right|du + \int_t^s\left| b_{\varepsilon}\left(u, \theta_{u, s}(y)\right) - b\left(u, \theta_{u, s}(y)\right)\right| du \leq C_{\theta} \sqrt{s - t}
	\label{equivalence of flows}
	\end{gather}
	by the Gr\"onwall inequality.
	
	\bigskip

For a given $\varepsilon > 0$, we say that the assumption \textbf{(A)} holds when conditions \textbf{(A1)} - \textbf{(A3)} are in force. Let us now introduce, under \textbf{(A)}, the quantities that will bound the difference of the
densities $p$ and $p_{\varepsilon}$ in our main result below. Set for $\varepsilon > 0$:

\begin{gather}\label{perturbation}
\Delta_{\varepsilon, b}(t,s) = \big\lVert\left|\left(b - b_{\varepsilon}\right)(u, \theta_{u,s}(y))\right|_{1}\big\rVert_{L_{1}} ,\\
\Delta_{\varepsilon, \sigma}(t,s) = \big\lVert\left|\left(\sigma - \sigma_{\varepsilon}\right)(u, \theta_{u,s}(y))\right|_{\gamma}\big\rVert_{L_{1}}, \\
\Delta_\varepsilon := \Delta_{\varepsilon, b} + \Delta_{\varepsilon, \sigma}.
\end{gather}
\bigskip
	Here $|\cdot|_{\beta}$ stands for the expression consisting of two terms:
	$$|f(t,x)|_{\beta}:=|f(t, x)|+\sup _{y \in\mathbb{R}^d, y \neq x} \frac{|f(t,x)-f(t,y)|}{|x-y|^{\beta}}.$$
	
	\bigskip
	
\end{itemize}

Now, for an arbitrary positive constant $\delta \in (\gamma/2, \gamma)$ and parameter $\alpha > 0$, we consider diagonal and off-diagonal maxima: $$ M_{\varepsilon, \alpha, \delta} = \max\limits_{\left(s - t\right) \leq \alpha} (s-t)^{\delta - \gamma/2}\left(\Delta_{\varepsilon}(t,s)\right)^{\gamma - \delta},$$
$$\bar{M}_{\varepsilon, \alpha, \delta} = \max\limits_{\left(s - t\right) \leq \alpha} \left(\Delta_{\varepsilon}(t,s)\right)^{\gamma - \delta},$$ 
$$M^{\mathcal{C}}_{\varepsilon, \alpha, \delta} = \max\limits_{\left(s - t\right) > \alpha} (s-t)^{\delta - \gamma/2}\left(\Delta_{\varepsilon}(t,s)\right)^{\gamma - \delta},$$
$$\bar{M}^{\mathcal{C}}_{\varepsilon, \alpha, \delta} = \max\limits_{\left(s - t\right) > \alpha} \left(\Delta_{\varepsilon}(t,s)\right)^{\gamma - \delta},$$

\bigskip

Let us state the main result of the paper.

\begin{theorem}\label{Main_thm}
	Let $p(t,s,x,y)$ and $p_{\varepsilon}(t,s,x,y)$ be the transition densities of \eqref{SDE_0} and \eqref{SDE_1} under \textbf{(A)}. Then there exists a constant $C > 0$ depending only on $T$ and parameters from the assumptions such that $\forall (t, x), (s, y) \in [0, T] \times \mathbb{R}^d$,
	\begin{gather}\label{Main_Est_L1}
		\left\lVert p-p_{\varepsilon}\right\rVert_{L_1^{\mu}\left(\mathbb{R}^d, L_1(\mathbb{R}^d)\right)}(t, s, x, y) \leq \dfrac{C}{\delta - \gamma/2} \cdot \left(M_{\varepsilon, \alpha, \delta} + M^{\mathcal{C}}_{\varepsilon, \alpha, \delta}\right) \\ \leq \dfrac{C}{\delta - \gamma/2}\left(\alpha^{\delta-\gamma/2}\bar{M}_{\varepsilon, \alpha, \delta} + T^{\delta - \gamma/2}\bar{M}^{\mathcal{C}}_{\varepsilon, \alpha, \delta}\right).
	\end{gather}	
\end{theorem}

\bigskip

\begin{remark}
	Suppose that $\Delta_{\varepsilon} \underset{\varepsilon \rightarrow 0}{\longrightarrow} 0$. Then the right-hand side of the last inequality in \eqref{Main_Est_L1} tends to zero for a suitable choice of $\alpha(\varepsilon) \rightarrow 0$. 
\end{remark}

\bigskip

In the following we will denote by $\langle \cdot,\cdot\rangle$ and $|\cdot| $ the Euclidean scalar product and the norm on $\R^d$. Also we keep the notation $D^\gamma_x = \prod\limits_{i = 1}^d D^{\gamma_i}_{x_i}$ to indicate the differentiation with respect to the multi-index $\gamma=(\gamma_1,\cdots,\gamma_d)\in \N^d $ and for which we denote  $|\gamma|=\sum_{i=1}^d \gamma_i$.

\section{McKean-Singer parametrix method} \label{McKean-Singer}

\subsection{Parametrix Representation of transition Densities}
Assume that \textbf{(A1)} and \textbf{(A2)} are in force.  These assumptions imply the existence of transition densities of \eqref{SDE_0} and \eqref{SDE_1} and, therefore,  will allow us to apply the PDE technique, namely, find the transition densities of processes solving \eqref{SDE_0} and \eqref{SDE_1} as fundamental solutions of the corresponding Kolmogorov equations.

\bigskip

Now let $\tilde{X}_{u, t}^{(\tau, z)}(x)$ denote the process starting at $x$ at time $t$ with dynamics
\begin{equation}\label{SDE_2}
d \tilde{X}_{u}^{(\tau, z)}=b\left(u, \theta_{u, \tau}(z)\right) du + \sigma\left(u, \theta_{u, \tau}(z)\right) d W_{u}, \ u \in [t, T].
\end{equation}

Usually there are two different choices for parameters $(\tau, z)$. The first possibility is to identify $(\tau, z)$ with the \textquotesingle \textquotesingle frozen\textquotesingle \textquotesingle \ initial point $(t,x)$ of the process \eqref{SDE_2}. However, we will choose the second possibility and freeze the terminal point $y$ at time $s$, letting $(\tau, z) =  (s, y)$. This approach is called \textit{backward} parametrix. 
\bigskip

In fact, after freezing the parameters $(\tau, z) =  (s, y)$,  \eqref{SDE_2} is inhomogeneous Gaussian process with mean and covariance matrix given ,respectively, by
\begin{equation}
	\tilde{\theta}_{s,t}(x) = x + \int\limits_{t}^{s} b(u, \theta_{u, s}(y))du = x + y - \theta_{t, s}(y),
\end{equation}
\begin{equation}
	\tilde{\mathcal{C}}_{s,t}(x)=\int\limits_t^s\sigma\sigma^*(u,\theta_{u,s}(y))du. 
\end{equation}

The idea is to use the density of the Gaussian process in \eqref{SDE_2}
to derive an estimate on the transition density $p^{\zeta}(t,s,x,y)$ of the process with  mollified versions $b_{\zeta}(t,x)$ and $\sigma_{\zeta}(t,x)$ of $b(t,x)$ and $\sigma(t,x)$. After that, the parametrix expansion for the transition density of the process \eqref{SDE_0} can be obtained by taking $\zeta \rightarrow 0$ \cite{DM, MPZ}. In this case, $\tilde{p}(t,s,x,y)$ is called a \textit{proxy}.
\bigskip

An important property is that the Gaussian transition density $\tilde{p}(t,s,x,y)$ of \eqref{SDE_2} satisfies the backward Kolmogorov equation:

\begin{equation}\label{backward}
	\left\{\begin{array}{l}
		\partial_{u} \tilde{p}(u, s, x, y)+\tilde{L}_{u} \tilde{p}(u, s, x, y)=0, \ t \leq u < s, \ x,y \in \mathbb{R}^{d} \\
		\tilde{p}(u, s, ., y) \underset{u \uparrow s}{\longrightarrow} \delta_{y}(.),	
\end{array}\right.
\end{equation}
where $$\tilde{L}_{u} \varphi(z, y)=\frac{1}{2} \operatorname{Tr}\left(\sigma \sigma^{*}(u, \theta_{u, s}(y)) D_{z}^{2} \varphi(z, y)\right)+\left\langle b(u, \theta_{u, s}(y)), D_{z} \varphi(z, y)\right\rangle$$
denotes the infinitesimal generator of \eqref{SDE_2} with frozen $(s, y) \in [0, T] \times \mathbb{R}^d$.
\bigskip

We now introduce a mollification procedure for the coefficients in \eqref{SDE_0}, i.e., setting $b_{\zeta} = b \star \xi_{\zeta}, \ \sigma_{\zeta} = \sigma \star \xi_{\zeta}$, where $\xi_{\zeta}(\cdot)=\zeta^{-(d+1)}\xi (\cdot/\zeta) $, $\xi\in C_0^{\infty}(\R^{d+1},\R_+),\ \int_{\R^{d+1}}\xi(z)dz=1 $, is a time-space mollifier  on $\R^{d+1} $ and $\star$ denotes the time-space convolution. It then follows from the H\"ormander theorem that the equation with dynamics
\begin{equation}
d X_{t}^{\zeta}=b_{\zeta}\left(t, X^{\zeta}_{t}\right) d t+\sigma_{\zeta}\left(t, X^{\zeta}_{t}\right) d W_{t}, \quad t \in[0, T],
\end{equation}
admits a smooth transition density $p^{\zeta}(t,s,x,y)$ for $t < s$. Therefore, $p^{\zeta}(t,s,x,y)$ must satisfy the forward Kolmogorov equation:

\begin{equation}\label{forward}
	\left\{\begin{array}{l}
		\partial_{u} p^{\zeta}(t, u, x, z)-L_{\zeta, u}^{*} p^{\zeta}(t, u, x, z)=0, \ t < u \leq s,\ z \in \mathbb{R}^{d} \\
		p^{\zeta}(t, u, x, .) \underset{u \downarrow t}{\longrightarrow} \delta_{x}(.),
	\end{array}\right.
\end{equation}
where $L_{\zeta}^*$ is an adjoint operator for the infinitesimal generator \eqref{generator}.

\bigskip

Our aim in this step is to estimate the transition density $p^\zeta(t,s,x,y)$ at every point using the Gaussian proxy \eqref{SDE_2}. To this end, we introduce a convolution, which will play a crucial role in our analysis:

\begin{equation}\label{conv_def}
	f \otimes g(t, s, x, y)=\int_{t}^{s} d u \int_{\mathbb{R}^{d}} d z f(t, u, x, z) g(u, s, z, y).
\end{equation}

Equations \eqref{backward}, \eqref{forward} yield the formal expansion, which is initially due to McKean and Singer (\cite{McS}):

$$\begin{aligned}\label{decomposition}
	\left(p^{\zeta}-\tilde{p}\right)(t, s, x, y)  &=\int_{t}^{s} d u \partial_{u} \int_{\mathbb{R}^{d}} d z p^{\zeta}(t, u, x, z) \tilde{p}(u, s, z, y) 
	 \\ &=\int_{t}^{s} d u \int_{\mathbb{R}^{d}} dz\left(\partial_{u} p^{\zeta}(t, u, x, z) \tilde{p}(u, s, z, y)+p^{\zeta}(t, u, x, z) \partial_{u} \tilde{p}(u, s, z, y)\right)  
	\\ &=\int_{t}^{s} d u \int_{\mathbb{R}^{d}} d z\left(L_{\zeta, u}^{*} p^{\zeta}(t, u, x, z) \tilde{p}(u, s, z, y)-p^{\zeta}(t, u, x, z) \tilde{L}_{u} \tilde{p}(u, s, z, y)\right)  
	\\ &=\int_{t}^{s} d u \int_{\mathbb{R}^{d}} d z p^{\zeta}(t, u, x, z)\left(L_{\zeta, u}-\tilde{L}_{u}\right) \tilde{p}(u, s, z, y) \\
	&:= p^{\zeta} \otimes H^{\zeta}(t, s, x, y),
\end{aligned}$$
where we used the Newton-Leibniz formula and the Dirac convergence for the first equality, equations \eqref{backward} and \eqref{forward} for the third one. We finally take the adjoint for the fourth equality. The function $H^{\zeta}(t,s,x,y) = \left(L_{\zeta, t}-\tilde{L}_{t}\right) \tilde{p}(t, s, x, y)$ is called a \textit{parametrix kernel}. With these notations the equality above rewrites as \begin{gather}
	p^{\zeta}(t,s,x,y) = \tilde{p}(t,s,x,y) + p^{\zeta} \otimes H^{\zeta}(t, s, x, y) = \tilde{p}(t,s,x,y) + \tilde{p} \otimes H^{\zeta}(t, s, x, y) + p^{\zeta} \otimes H^{\zeta, 2}(t, s, x, y) \implies \\
	p^{\zeta}(t,s,x,y) = \sum\limits_{i = 0}^{+\infty} \tilde{p} \otimes H^{\zeta, i}(t,s,x,y),
\end{gather} 

where $\tilde{p} \otimes H^{\zeta, i}(t,s,x,y) = \left[\tilde{p} \otimes H^{\zeta, i - 1}\right] \otimes H^{\zeta}(t,s,x,y)$ and $\tilde{p} \otimes H^{\zeta, 0}(t,s,x,y) := \tilde{p} (t,s,x,y).$

\bigskip
 Now it can be shown that $H^\zeta \underset{\zeta \rightarrow 0}{\rightarrow} H := \left(L_{t}-\tilde{L}_{t}\right) \tilde{p}(t, s, x, y) $ pointwise. Therefore, by the bounded convergence theorem, $p^\zeta(t,v,x,y)\underset{\zeta\rightarrow 0}{\longrightarrow} \tilde p(t,s,x,y)+\sum_{r = 1}^{\infty}\tilde p\otimes H^{r}(t,s,x,y)$. Moreover, Theorem 11.1.4 in \cite{stroock} implies that for any bounded continuous function $f:\R^d\rightarrow \R$,
 $$\mathbb E\left[f(X^{\zeta}_u)\right] \underset{\zeta \rightarrow 0}{\longrightarrow} \mathbb E\left[f(X_u)\right].$$
 This eventually gives  (see the details in \cite{KKM})
 \begin{equation}\label{par_series}
 p(t,s,x,y)=\tilde p(t,s,x,y)+\sum_{r = 1}^{\infty}\tilde p\otimes H^{r}(t,s,x,y).
 \end{equation}
Recall that the decomposition \eqref{par_series} is formal until we have proved the uniform convergence of the series on the right-hand side.

\subsection{Flow associated with drift}

In this subsection, we derive the properties of the flow solving the ODE \eqref{flow}. These properties will be important for further analysis.

\begin{proposition}[Lipschitz continuity of the flow]\label{Lipschitz_flow}
$\theta_{t, s}(\cdot)$ is Lipschitz continuous in the spatial variable, i.e., there exists $C > 0$ such that $\forall x, y \in \mathbb{R}^d$ 
$$\left|\theta_{t, s}(x) - \theta_{t, s}(y)\right| \leq C \left|x - y\right|.$$

\end{proposition}
\proof
 Let's fix the initial and terminal times $t$ and $s$, respectively. Let $x, y \in \mathbb{R}^d$. Then
\begin{gather}
	\left|\theta_{t, s}(x) - \theta_{t, s}(y)\right| = \left|x - y + \int_t^s \left(b(u, \theta_{u,s}(x)) - b(u, \theta_{u,s}(y))\right)du\right| \leq \\ \leq \left|x - y\right| + K \cdot \int\limits_{t}^{s}\left|\theta_{u, s}(y) - \theta_{u, s}(x)\right|du.
\end{gather}

By the Gr\"onwall inequality, 
\begin{equation}
	\left|\theta_{t, s}(x) - \theta_{t, s}(y)\right| \leq \left|x - y\right| \cdot \exp(C(s - t)) \leq \left|x - y\right| \cdot \exp(CT).
\end{equation}
\endproof
\bigskip

It is well-known, that under assumption \textbf{(A2)}, the flow $\theta_{t,s}(\cdot)$ has Lipschitz inverse and enjoys the semigroup property, i.e., 
$$\theta_{t,s}\left(\theta_{s,u}(y)\right) = \theta_{t,u}(y),$$
$$\theta_{s,t}\left(\theta_{t,u}(y)\right) = \theta_{s,u}(y).$$
When $t = u$ or $s = u$, respectively, the equalities above take the form 
\begin{gather}\label{Semigroup_property}
\theta_{t,s}\left(\theta_{s,t}(y)\right) = y,
\theta_{s,t}\left(\theta_{t,s}(y)\right) = y.
\end{gather}

The direct consequence of the semigroup property is the next proposition (see \cite{DM}):
\begin{proposition}[Bi-Lipschitz property of flow]
	The deterministic flow $\theta_{u, s}(\cdot)$ solving \eqref{flow} enjoys the bi-Lipschitz property, i.e., there exists a constant $C \geq 1$ such that for all $x, y \in \mathbb{R}^{d}$ 
	\begin{equation}
		C^{-1}\left|y - \theta_{s, u}(x)\right| \leq \left|x - \theta_{u, s}(y)\right| \leq C\left|y - \theta_{s, u}(x)\right|.
	\end{equation}
\end{proposition}

\bigskip

\subsection{Convergence of parametrix series for diffusions with unbounded drift}

Recall that the transition density of the process \eqref{SDE_0} can be expanded into the
formal parametrix series \eqref{par_series}. We recall for the sake of completeness the key steps in the proof of its uniform convergence. In this subsection, we suppose that the coefficients of \eqref{SDE_0} satisfy \textbf{(A1)}, \textbf{(A2)}.

\bigskip
From direct computations, for all $0 \leq t < s \leq T, \ x,y \in \mathbb{R}^d$ and any multiindex $\nu$ with $|\nu| \leq 4 $, there exists $C > 0$ such that 
\begin{equation}
	\left|D^{\nu}_x\tilde p(t,s,x,y)\right| \leq \dfrac{C}{\left(s-t\right)^{(|\nu| + d)/2}} \cdot \exp\left(\frac{-\displaystyle\left|\theta_{t,s}(y)-x\right|^2}{\displaystyle C(s-t)}\right).
	\label{est_der_dens}
\end{equation}

Next, applying H\"older continuity of the coefficients and \eqref{est_der_dens}, we readily get that there is a constant $C > 0$ such that
\begin{equation}
	\displaystyle |H(t,s,x,y)| \leq \frac{C}{(s-t)^{d/2+1-\gamma/2}} \cdot \exp\left(\frac{-\displaystyle\left|\theta_{t,s}(y)-x\right|^2}{\displaystyle C(s-t)}\right).
	\label{est_kernel}
\end{equation} 
\bigskip

The next step is estimating the convolutions of higher order. Direct iteration of \eqref{decomposition} leads to unbounded growth of constant $C$ and, as a consequence, \eqref{par_series} will diverge. The following proposition allows one to manage the iteration procedure without the constants deteriorating  at every step. This estimate was told us by S.Menozzi (personal communication).  
For the readers convenience we give a sketch of the proof. 

\begin{proposition}\label{majorizing} For any $C$, $0 < C < \infty$
	the expression $$\frac{C}{(s-t)^{d/2}} \cdot \exp\left(\frac{-\displaystyle\left|\theta_{t,s}(y)-x\right|^2}{\displaystyle C(s-t)}\right)$$ is a lower bound for the transition density $\bar{p}$ of the auxiliary diffusion process $\bar{X}_t$ with the dynamics
	\begin{equation}
		d\bar{X}_t = b\left(t, \bar{X}_t\right) d t+\lambda I d W_{t},\ t \in[0, T],
		\label{Maj_dens} 
	\end{equation}
	where $\lambda = \lambda(C)$ is a positive large parameter depending on the quantities in the assumptions. In other words, the transition density $\bar{p}$ of the process $\bar{X}_t$ satisfies $$
	\bar{C} \bar{p}(t, s, x, y) \geq C(s-t)^{-\frac{d}{2}} \exp \left(-\frac{\left|\theta_{s, t}(x)-y\right|^{2}}{C(s-t)}\right)
	$$ for a constant $\bar{C}$.
	
	\bigskip
	Importantly, $\bar p $ satisfies the Chapman-Kolmogorov identity. Namely,
	\begin{equation}
		\int_{\R^d} \bar p(t,v,x,z)\bar p(v,s,z,y)dz=\bar p(t,s,x,y).
	\end{equation}
	It also holds that
	\begin{equation}\label{UP_BD_BAR_P}
		\bar{p}(t,s,x,y)\le \frac{C}{(s-t)^{\frac{d}{2}}}\exp\left(-\frac{|\theta_{t,s}(y)-x|^2}{C(s-t)}\right),
	\end{equation}
	for the flow $\theta_{t,s}(y) $ chosen above that solves \eqref{flow}, up to a possible modification of $C$.	
\end{proposition}
\proof
Let us consider 
\begin{equation}
	C_{T}(s-t)^{-\frac{d}{2}} \exp \left(-\frac{\left|\theta_{s, t}(x)-y\right|^{2}}{C_{T}(s-t)}\right)	
	\label{1}
\end{equation} 
with $C_{T}$ depending only on \textbf{(A)}. We have to show that this expression can serve as a lower bound for an auxiliary diffusion $\bar{X}_{t}$ with sufficiently large $\lambda\left(C_{T}\right)$. For the linear drift $b(t, x)=G_{t} x$, it follows directly from (3.2), (3.4) and (3.7) in Proposition 3.1 in \cite{DM}. Indeed, in this case the Gaussian density $\bar{p}(t, s, x, y)$ can be written explicitly
$$
\bar{p}(t, s, x, y)=(2 \pi)^{-\frac{d}{2}} \operatorname{det}^{-\frac{1}{2}}\left(K_{s}\right) \exp \left(-\frac{1}{2}\left\langle K_{s}^{-1}(R(s, t) x-y), R(s, t) x-y\right\rangle\right),
$$
where
$$
K_{s}=\lambda^{2} \int_{t}^{s} R(s . u) R^{*}(s, u) d u:=\lambda^{2} \bar{R}(s, t), \quad K_{s}^{-1}=\lambda^{-2} \bar{R}^{-1}(s, t).
$$
$R(s, t)_{0 \leq t, s \leq T}$ stands for the resolvent associated with $\left(G_{t}\right)_{0 \leq t \leq T}, \ R(s, t) x=\theta_{s, t}(x) .$ The density $\bar{p}(t, s, x, y)$ can be rewritten as
$$
\bar{p}(t, s, x, y)=(2 \pi)^{-\frac{d}{2}} \lambda^{-d} \operatorname{det}^{-1 / 2} \bar{R}(s, t) \exp \left(-\frac{1}{2 \lambda^{2}}\left\langle\bar{R}^{-1}(s, t)\left(\theta_{s, t}(x)-y\right), \theta_{s, t}(x)-y\right\rangle\right) \geq
$$ \begin{equation}
K_{T}(s-t)^{-\frac{d}{2}} \lambda^{-d} \exp \left(-\frac{\left|\theta_{s, t}(x)-y\right|^{2}}{\lambda^{2} K_{T}(s-t)}\right).
\label{2}
\end{equation} Now we take $\lambda$ such that $\lambda^{2} K_{T}>C_{T}$ and then with this fixed $\lambda$ take $\bar{C}$ such that $\bar{C} K_{T} \lambda^{-d} > C_{T}$. Then we obtain from \eqref{1} and \eqref{2}
$$\bar{C} \bar{p}(t, s, x, y) \geq \bar{C} K_{T}(s-t)^{-\frac{d}{2}} \lambda^{-d} \exp \left(-\frac{\left|\theta_{s, t}(x)-y\right|^{2}}{\lambda^{2} K_{T}(s-t)}\right) \geq $$ $$
	C_{T}(s-t)^{-\frac{d}{2}} \exp \left(-\frac{\left|\theta_{s, t}(x)-y\right|^{2}}{C_{T}(s-t)}\right). $$
The general case uses linearization and follows the strategy developed in \cite{DM}, pages $1599-1613 .$ Many objects become simple or trivial because of a constant diffusion matrix $\lambda I_{d}$. We only indicate the key points resulting to a desire majorization. Inequality of Proposition 4.2 becomes
$$
\sup _{0 \leq s \leq t}\left|\varphi_{s}\right|^{2} \leq \frac{C}{\lambda} \cdot \frac{\left|y-\theta_{t}(x)\right|^{2}}{t}
$$ with a constant $C$ depending on \textbf{(A)} and $T$ only. The control of Proposition 4.3 has the following form
$$
\mathbb P\left\{\forall t \in [0, T-\varepsilon],\ \begin{aligned}
	\left|\Theta_{t}-\Gamma_{t}\right| & \leq C(\mu)(T-t)^{-\frac{1}{2}+\frac{\eta}{8}} \\
	\left|v_{t}^{0}-\gamma_{t}\right| & \leq \frac{C(\mu)}{\lambda^{2}}(T-t)^{-\frac{1}{2}+\frac{\eta}{8}}
\end{aligned}\right\} \geq 1-\mu, 
$$
and (4.27) and (4.29) become $$ \mathbb{E}\left(\mathbb{I}_{\bar{C}} \sum_{i=1}^{6} R_{T-\varepsilon}^{i}\right) \leq C_{1}(\lambda)+\frac{C_{2}}{\lambda} \frac{\left|\theta_{T}\left(x_{0}\right)-y_{0}\right|^{2}}{T} $$
$$	-\ln \bar{p}(t, s, x, y) \leq \frac{d}{2} \ln (s-t)+C_{1}(\lambda)+\frac{C_{2}}{\lambda} \frac{\left|\theta_{T}\left(x_{0}\right)-y_{0}\right|^{2}}{T},
$$ where $C_2$ is a constant depending only on the quantities in the assumptions. Here the parameter $\varepsilon$ has the same meaning as in \cite{DM}, not to be confused with the parameter that we used for the perturbed equation \eqref{SDE_1}.

Hence, for a constant $C_{T}$ in the upper bound for $p(t, s, x, y)$ we can choose sufficiently large $\lambda$ and then sufficiently large $\bar{C}$ such that
$$\bar{C} \bar{p}(t, s, x, y) \geq \bar{C} \exp \left(-C_{1}(\lambda)\right)(s-t)^{-\frac{d}{2}} \exp \left(-\frac{C_{2}}{\lambda} \frac{\left|\theta_{s, t}(x)-y\right|^{2}}{s-t}\right) \geq $$
$$	C_{T}(s-t)^{-\frac{d}{2}} \exp \left(-\frac{\left|\theta_{s, t}(x)-y\right|^{2}}{C_{T}(s-t)}\right). $$
\endproof

\begin{corollary}
	By Theorem 1.1 in \cite{DM} $$
	p(t, s, x, y) \leq C_{T}(s-t)^{-\frac{d}{2}} \exp \left(-\frac{\left|\theta_{s, t}(x)-y\right|^{2}}{C_{T}(s-t)}\right)
	$$
	with $C_{T}$ depending only on \textbf{(A)}. Proposition \ref{majorizing} implies
	$$
	\bar{C} \bar{p}(t, s, x, y) \geq p(t, s, x, y)
	$$
	for sufficiently large constant $\bar{C}$.
\end{corollary}
\bigskip

\begin{remarks} ~\\

	\begin{itemize}

\item[\textbf{1.}] In fact, these modifications of the proofs allow us to construct a majorizing diffusion for a broader class of degenerate Kolmogorov-type systems considered in \cite{DM}.

\item[\textbf{2.}] It should also be noted that in \cite{AMP}, when studying Geometric Average Asian options, the authors considered Kolmogorov-type operators and obtained two-sided estimates for the fundamental solution, 
these estimates are also given in terms of the fundamental solutions of the corresponding equations with constant  diffusions. 
	
\item[\textbf{3.}]
Sometimes, in our paper, we will automatically renew the majorizing density $\bar{p}(t,s,x,y)$, meaning the following sequence of inequalities:
$$\left|\theta_{t,s}(y) - x\right|^{\beta} \bar{p}(t,s,x,y) \leq \left|\theta_{t,s}(y) - x\right|^{\beta}\frac{C_1}{(s-t)^{d/2}}\exp\left(-\frac{|\theta_{t,s}(y)-x|^2}{C_1(s-t)}\right) \leq $$ $$\leq \frac{C_2}{(s-t)^{d/2}}\exp\left(-\frac{|\theta_{t,s}(y)-x|^2}{C_2(s-t)}\right) \leq \bar{p}(t,s,x,y).$$

\end{itemize}
\end{remarks}
We can now rewrite the estimates \eqref{est_der_dens} and \eqref{est_kernel} in the appropriate form using the majorizing density $\bar{p}(t,s,x,y)$:
\begin{equation}
	\left|\tilde p(t,s,x,y)\right| \leq C \cdot \bar{p}(t,s,x,y), 
	\label{est_den_new} 
\end{equation}
\begin{equation}
	\displaystyle |H(t,s,x,y)| \leq \dfrac{C}{(s-t)^{1-\gamma/2}} \cdot \bar{p}(t,s,x,y).
	\label{est_ker_new} 
\end{equation}

The following proposition shows that the same upper bounds can be derived for the transition density $\tilde{p}_{\varepsilon}(t,s,x,y)$ if we additionally assume that \textbf{(A3)} holds.

\begin{proposition}
	The transition density $\tilde{p}_{\varepsilon}(t,s,x,y)$ of the perturbed process \eqref{SDE_1} admits an upper bound  \eqref{est_den_new} with the same majorizing density $\bar{p}$:
		$$\left|\tilde p_{\varepsilon}(t,s,x,y)\right| \leq C \cdot \bar{p}(t,s,x,y). $$
\end{proposition}
\proof
Similarly to \eqref{est_der_dens}, $\tilde{p}_{\varepsilon}(t,s,x,y)$ satisfies the following estimate:
$$
\left|\tilde{p}_{\varepsilon}(t, s, x, y)\right| \leq \frac{C}{(s-t)^{d / 2}} \cdot \exp \left(\frac{-\left|\theta_{t, s}^{\varepsilon}(y)-x\right|^{2}}{C(s-t)}\right)=\frac{C}{(s-t)^{d / 2}} \cdot \exp \left(\frac{-\left|\theta_{t, s}^{\varepsilon}(y)-\theta_{t, s}(y)+\theta_{t, s}(y)-x\right|^{2}}{C(s-t)}\right) \leq 
$$
$$\leq \frac{2C}{(s-t)^{d / 2}} \cdot \exp \left(\frac{-\left|\theta_{t, s}(y)-x\right|^{2}}{2C(s-t)}\right) \exp\left(\frac{\left|\theta_{t, s}^{\varepsilon}(y)-\theta_{t, s}(y)\right|^{2}}{C(s-t)}\right).$$

The factor $\exp\left(\frac{\left|\theta_{t, s}^{\varepsilon}(y)-\theta_{t, s}(y)\right|^{2}}{C(s-t)}\right)$ is bounded due to the consequence of the assumption \textbf{(A3)}, while the remaining factor can be viewed as a lower bound for majorizing density.  
\endproof
\bigskip

Recall that the proposition above immediately implies the control of the perturbed parametrix kernel $H_{\varepsilon}$:
$$\left|H_{\varepsilon}(t,s,x,y)\right| \leq \dfrac{C}{(s-t)^{1 - \gamma/2}}\bar{p}(t,s,x,y).$$
\bigskip

Using the Chapman-Kolmogorov equation, we get that there is an upper bound for the $n$-order convolution of the parametrix kernels with the transition densities $\tilde{p}(t,s,x,y)$ and $\tilde{p}_{\varepsilon}(t,s,x,y)$:
	\begin{gather}
	\left|\tilde{p} \otimes H^{n}\right|(t,s,x,y)	\leq C^{n + 1} \dfrac{\Gamma\left(\frac{\gamma}{2}\right)^{n}}{\Gamma(1 + \frac{n\gamma}{2})} (s-t)^{n\gamma/2} \cdot \bar{p}(t, s, x, y), \\
	\left|\tilde{p}_{\varepsilon} \otimes H^{n}_{\varepsilon}\right|(t,s,x,y)	\leq C^{n + 1} \dfrac{\Gamma\left(\frac{\gamma}{2}\right)^{n}}{\Gamma(1 + \frac{n\gamma}{2})} (s-t)^{n\gamma/2} \cdot \bar{p}(t, s, x, y).
	\label{conv_n-order_est}
\end{gather}

\bigskip

The above estimate \eqref{conv_n-order_est} means that parametrix series \eqref{par_series} converges absolutely and uniformly on $\left([0, T] \times \mathbb{R}^d\right)^2$: 
\begin{equation}
	\sum\limits_{n=0}^{\infty}\left|\tilde{p} \otimes H^{n}\right|(t,s,x,y) \leq \frac{C}{(s-t)^{1-\gamma/2}} \cdot \bar{p}(t,s,x,y)
\end{equation}
for a constant $C > 0$. The same argument is valid for the parametrix series of $p_{\varepsilon}$.
\bigskip

\section{Stability of parametrix series}
\subsection{$L_1$ - perturbation case}

In this subsection, we will investigate the stability of \eqref{SDE_0} under $L_{1}$- type perturbation \eqref{perturbation} through the difference of the respective series. For a given fixed parameter $\varepsilon$, under  $\textbf{(A1)} - \textbf{(A2)}$, the densities $p(t, s, x, y), p_{\varepsilon}(t, s, x, y)$ of the processes in \eqref{SDE_0}, \eqref{SDE_1} both admit a parametrix expansion of the type \eqref{par_series}. This  decomposition allows us to get an upper bound for the difference of the transition densities using \textquotesingle \textquotesingle term-by-term\textquotesingle \textquotesingle \ estimation. 
\bigskip

Clearly, 

$$\begin{aligned}\label{diff_series}
 \left|p(t,s,x,y) - p_{\varepsilon}(t,s,x,y)\right| & = \left|\sum_{r=0}^{\infty} \tilde{p} \otimes H^{r}(t,s,x,y) - \tilde{p}_{\varepsilon} \otimes H^{r}_{\varepsilon}(t,s,x,y)\right| \\ & \leq \left|\tilde{p}(t,s,x,y) - \tilde{p}_{\varepsilon}(t,s,x,y)\right| + \sum_{r=1}^{\infty} \left|\tilde{p} \otimes H^{r}(t,s,x,y) - \tilde{p}_{\varepsilon} \otimes H^{r}_{\varepsilon}(t,s,x,y)\right|.
\end{aligned}$$

The strategy is to study the above difference by using some well-known properties of the Gaussian kernels and their derivatives. First, we study the difference of the \textit{main terms}. From now on, we fix a positive parameter $\delta \in (\gamma/2, \gamma)$, which will be necessary  to avoid singualrities when estimating the differences of the higher order convolutions. 
\bigskip
Let us introduce some notations:

$$\phi(t,s; y) = \left(\int_{t}^{s}\left|\left(b-b_{\varepsilon}\right)\left(u, \theta_{u, s}(y)\right)\right|_{1} d u\right)^{\gamma-\delta} + \left(\int_{t}^{s}\left|\left(\sigma-\sigma_{\varepsilon}\right)\left(u, \theta_{u, s}(y)\right)\right|_{\gamma} d u\right)^{\gamma-\delta},
$$
\bigskip
$$\psi(t,s; y) = \left|\left(b-b_{\varepsilon}\right)\left(t, \theta_{t, s}(y)\right)\right|_{1} + \left|\left(\sigma-\sigma_{\varepsilon}\right)\left(t, \theta_{t, s}(y)\right)\right|_{\gamma}, $$

$$||f(t,s,x, y)||_{L_1\left(\mathbb{R}^d\right)} = \int_{\mathbb{R}^d}|f(t,s,x,y)| dy.$$
\bigskip

\begin{lemma}[Difference of the main terms and their derivatives]\label{diff_main_terms}
	There exists $C>0$ such that for all $0 \leq t < s \leq T$, $x, y \in \mathbb{R}^d$, $\delta \in (\gamma/2, \gamma)$, and $\nu, |\nu| \leq 4$:
	\begin{gather}
	\left|D_x^{\nu}\tilde{p}(t,s,x,y) - D_x^{\nu}\tilde{p}_{\varepsilon}(t,s,x,y)\right| \leq \\ \leq \dfrac{C\bar{p}(t,s,x,y)}{(s - t)^{|\nu|/2}} \left(\dfrac{\left(\int_{t}^{s}\left|\left(b-b_{\varepsilon}\right)\left(u, \theta_{u, s}(y)\right)\right| d u\right)^{\gamma-\delta}}{(s-t)^{\gamma/2 - \delta/2}} + \dfrac{\left(\int_{t}^{s}\left|\left(\sigma-\sigma_{\varepsilon}\right)\left(u, \theta_{u, s}(y)\right)\right| d u\right)^{\gamma-\delta}}{(s-t)^{\gamma - \delta}}\right).
	\end{gather}

	Consequently, in terms of $L_1$ norm, the difference of the main terms can be rewritten as  
	$$\begin{aligned} \label{L_p_diff_main_terms}
	&||\left(\tilde{p}_{\varepsilon} - \tilde{p}\right)(t,s,x,y)||_{L_1(\mathbb{R}^d)} \leq   C(s-t)^{(\gamma - \delta)/2} \left[\left(\Delta_{\varepsilon, b}^{(\gamma-\delta)}	\right)^{\gamma-\delta} + \left(\Delta_{\varepsilon, \sigma}^{(\gamma-\delta)}	\right)^{\gamma-\delta}\right]\end{aligned}.$$
	
\end{lemma}
\bigskip
The previous lemma estimates the difference of the main terms of the expansions. The next lemma quantifies
the difference between the parametrix kernels $H(t,s,x,y)$ and $H_{\varepsilon}(t,s,x,y)$.
\bigskip
\begin{lemma}[Difference of parametrix kernels]\label{diff_par_ker}
	For all $0 \leq t < s \leq T$ and $x, y \in \mathbb{R}^d$, there exists $C>0$ such that
	$$|H - H_{\varepsilon}|(t,s,x,y) \leq C\bar{p}(t,s,x,y) \left( \frac{\psi(t,s;y)}{(s-t)^{1-\gamma/2}}  + \frac{\phi(t,s;y)}{(s-t)^{1 + \gamma/2 - \delta}}\right).$$
\end{lemma}
\bigskip

For technical reasons, we postpone the proofs of Lemmas \ref{diff_main_terms} and \ref{diff_par_ker} till the Appendix. 
\bigskip

To complete the proof of the main theorem, we need to handle the difference of the $n$-th order convolutions. To achieve this, we split the  convolution into two parts:

\begin{gather}\label{decomp}
	\left|\tilde p \otimes H^{n+1} - \tilde p_{\varepsilon} \otimes H_{\varepsilon}^{n+1}\right|(t,s,x,y) \leq \\ \leq \left|\left(\tilde p_{\varepsilon} \otimes H_{\varepsilon}^{n}\right) \otimes \left(H - H_{\varepsilon}\right)\right|(t,s,x,y) + \left|\left(\tilde p \otimes H^{n} - \tilde p_{\varepsilon} \otimes H_{\varepsilon}^{n}\right) \otimes H\right|(t,s,x,y). 
\end{gather}

\bigskip

The next lemma allows one to control the first term of \eqref{decomp}.

\begin{lemma}\label{second_part}
For all $0 \leq t < s \leq T$ and $x, y \in \mathbb{R}^d$, there exists $C>0$ such that
$$\left|\left(\tilde p_{\varepsilon} \otimes H_{\varepsilon}^{n}\right) \otimes \left(H - H_{\varepsilon}\right)\right|(t,s,x,y) \leq C^{n + 2} (s-t)^{n\gamma/2}\dfrac{\Gamma^n(\gamma/2)}{\Gamma(1 + n\gamma/2)}\bar{p}(t,s,x,y)\int\limits_t^s \dfrac{\psi(u,s;y)}{(s-u)^{1 - \gamma/2}}du + $$$$+ C^{n + 2} (s-t)^{\delta + (n-1)\gamma/2}\dfrac{\Gamma^n(\gamma/2)\Gamma(\delta - \gamma/2)}{\Gamma(1 + \delta + (n - 1)\gamma/2)}\bar{p}(t,s,x,y)\phi(t,s;y).$$
\end{lemma}
\proof

Using the control of $n$-th order convolution \eqref{conv_n-order_est} and Lemma \ref{diff_par_ker}, we obtain:

$$\left|\left(\tilde p_{\varepsilon} \otimes H_{\varepsilon}^{n}\right) \otimes \left(H - H_{\varepsilon}\right)\right|(t,s,x,y) \leq $$ $$\leq C^{n + 2} \dfrac{\Gamma^n(\gamma/2)}{\Gamma(1 + n\gamma/2)} \bar{p}(t,s,x,y) \left((s-t)^{n\gamma/2} \int\limits_t^s \dfrac{\psi(u,s;y)}{(s-u)^{1 - \gamma/2}}du + \phi(t,s;y) \int\limits_t^s \dfrac{(u-t)^{n\gamma/2}}{(s-u)^{1 + \gamma/2 - \delta}}du\right) = $$ $$= C^{n + 2} (s-t)^{n\gamma/2}\dfrac{\Gamma^n(\gamma/2)}{\Gamma(1 + n\gamma/2)}\bar{p}(t,s,x,y)\int\limits_t^s \dfrac{\psi(u,s;y)}{(s-u)^{1 - \gamma/2}}du + $$$$+ C^{n + 2} (s-t)^{\delta + (n-1)\gamma/2}\dfrac{\Gamma^n(\gamma/2)\Gamma(\delta - \gamma/2)}{\Gamma(1 + \delta + (n - 1)\gamma/2)}\bar{p}(t,s,x,y)\phi(t,s;y).$$

\endproof

\bigskip

To get an upper bound for the second term of \eqref{decomp}, we will iterate it. We estimate the difference between convolutions of the first order.    
\bigskip

\begin{lemma}\label{diff_first_ord_conv}
	For all $0 \leq t < s \leq T$ and $x, y \in \mathbb{R}^d$, there exists $C>0$ such that
	$$\left|\tilde p \otimes H - \tilde p_{\varepsilon} \otimes H_{\varepsilon}\right|(t,s,x,y) \leq \dfrac{C^2}{\delta - \gamma/2} \bar{p}(t,s,x,y)\left(\dfrac{B\left(1 + \delta - \gamma; \gamma/2\right)}{(s-t)^{\gamma/2 - \delta}} \phi(t,s;y) + \int\limits_t^s \dfrac{\psi(u,s;y)}{(s-u)^{1 - \gamma/2}}du \right).$$
\end{lemma}
\proof

Using the strategy described above, we write
\
\begin{equation}
	\left|\left(\tilde p \otimes H - \tilde p_{\varepsilon} \otimes H_{\varepsilon}\right)(t,s,x,y)\right| \leq \left|\left(\tilde p - \tilde p_{\varepsilon}\right) \otimes H\right|(t,s,x,y) +\left| \tilde p_{\varepsilon} \otimes \left(H - H_{\varepsilon}\right)\right|(t,s,x,y).
\end{equation}
 \bigskip
 
 To estimate the first term in the expression above, using the semigroup property, we write
 
 $$\left|b - b_{\varepsilon}\right|(t, \theta_{t,u}(z)) \bar{p}(u,s,z,y)\leq \left[ \left|\left(b - b_{\varepsilon}\right)(t, \theta_{t,s}(\theta_{s,u}(z))) - \left(b - b_{\varepsilon}\right)(t, \theta_{t,s}(y))\right| + \right.$$ $$\left.+\left|\left(b - b_{\varepsilon}\right)(t, \theta_{t,s}(y))\right|\right]\bar{p}(u,s,z,y) \leq$$ $$\leq C\left( \left|\theta_{s,u}(z) - y\right| + 1\right)\left|\left(b - b_{\varepsilon}\right)(t, \theta_{t,s}(y))\right|_{1} \bar{p}(u,s,z,y) \leq C\left|\left(b - b_{\varepsilon}\right)(t, \theta_{t,s}(y))\right|_{1} \bar{p}(u,s,z,y).$$
 \bigskip

 Recall that the same argument is valid for the diffusion coefficients $\sigma$, $\sigma_{\varepsilon}$: $$\left|\sigma - \sigma_{\varepsilon}\right|(t, \theta_{t,u}(z)) \bar{p}(u,s,z,y)\leq C\left|\left(\sigma - \sigma_{\varepsilon}\right)(t, \theta_{t,s}(y))\right|_{\gamma} \bar{p}(u,s,z,y).$$
 \bigskip

 Now, from Lemma \ref{diff_main_terms} and \eqref{est_ker_new},
 
 $$\left|\left(\tilde p - \tilde p_{\varepsilon}\right) \otimes H\right|(t,s,x,y) \leq $$

 $$C^2 \bar{p}(t,s,x,y)\left[\left(\int_t^s\left|\left(b-b_{\varepsilon}\right)\left(u, \theta_{u, s}(y)\right)\right|_{1} d u\right)^{\gamma-\delta} \int_t^s (u - t)^{\delta/2 - \gamma/2}(s-u)^{\gamma/2 - 1}du + \right.$$ $$\left.+ \left(\int_{t}^{s}\left|\left(\sigma-\sigma_{\varepsilon}\right)\left(u, \theta_{u, s}(y)\right)\right|_{\gamma} d u\right)^{\gamma-\delta}\int^s_t (u - t)^{\delta - \gamma}(s-u)^{\gamma/2 - 1}du \right] \leq$$ $$\leq \dfrac{C^2B\left(1 + \delta - \gamma; \gamma/2\right)}{(s-t)^{\gamma/2 - \delta}} \bar{p}(t,s,x,y) \left[\left(\int_t^s\left|\left(b-b_{\varepsilon}\right)\left(u, \theta_{u, s}(y)\right)\right|_{1} d u\right)^{\gamma-\delta}  + \right.$$ $$\left. \left(\int_{t}^{s}\left|\left(\sigma-\sigma_{\varepsilon}\right)\left(u, \theta_{u, s}(y)\right)\right|_{\gamma} d u\right)^{\gamma-\delta} \right] =   \dfrac{C^2B\left(1 + \delta - \gamma; \gamma/2\right)}{(s-t)^{\gamma/2 - \delta}} \bar{p}(t,s,x,y)\phi(t,s;y).$$
 
 Combining the upper bound derived above with Lemma \ref{diff_par_ker} completes the proof.
 \endproof
 \bigskip 
 
 Now we are ready to prove the main theorem \ref{Main_thm} of this paper.
 
 \begin{proof}
 	From Lemmas \ref{diff_first_ord_conv}, \ref{second_part} and Jensen's inequality,  we have the following $L_1 - L_1$ upper bounds:
 	
$$||\left(\tilde p \otimes H - \tilde p_{\varepsilon} \otimes H_{\varepsilon}\right)(t,s,x,y)||_{L_1^{\mu}(\mathbb{R}^d, L_1(\mathbb{R}^d))} \leq$$ $$\leq \dfrac{C^2}{\delta-\gamma/2}B\left(1 + \delta - \gamma; \gamma/2\right)(s-t)^{(\gamma + \delta)/2}\left[\left(\Delta_{\varepsilon, b}	\right)^{\gamma-\delta} + \left(\Delta_{\varepsilon, \sigma}	\right)^{\gamma-\delta}\right],$$

$$||\left(\left(\tilde p_{\varepsilon} \otimes H_{\varepsilon}^{n}\right) \otimes \left(H - H_{\varepsilon}\right)\right)(t,s,x,y)||_{L_1^{\mu}(\mathbb{R}^d, L_1(\mathbb{R}^d))} \leq$$ $$\leq C^{n+2}(s-t)^{\delta + (n-1)\gamma/2}\dfrac{\Gamma^n(\gamma/2)}{\Gamma(1 + n\gamma/2)}\left[\left(\Delta_{\varepsilon, b}	\right)^{\gamma-\delta} + \left(\Delta_{\varepsilon, \sigma}	\right)^{\gamma-\delta}\right].$$

\bigskip

The next step is the iteration of  \eqref{decomp} using the above estimates:

$$||\left(\tilde p \otimes H^{n+1} - \tilde p_{\varepsilon} \otimes H_{\varepsilon}^{n+1}\right)(t,s,x,y)||_{L_1^{\mu}(\mathbb{R}^d, L_1(\mathbb{R}^d))} \leq $$ $$\leq C^{n+2}(s-t)^{\delta + (n-1)\gamma/2}\dfrac{\Gamma^n(\gamma/2)}{\Gamma(1 + n\gamma/2)}\left[\left(\Delta_{\varepsilon, b}	\right)^{\gamma-\delta} + \left(\Delta_{\varepsilon, \sigma}	\right)^{\gamma-\delta}\right] + $$ $$+C \int_{\mathbb{R}^d}\int_{\mathbb{R}^{d}} \int_{t}^{s} d u \int_{\mathbb{R}^{d}}\left|\tilde{p} \otimes H^{n}-\tilde{p}_{\varepsilon} \otimes H_{\varepsilon}^{n}\right|(t, u, x, z) \frac{\bar{p}(u, s, z, y)}{(s-u)^{1-\gamma/2}} d z d y \ \mu(dx) = $$ $$= C^{n+2}(s-t)^{\delta + (n-1)\gamma/2}\dfrac{\Gamma^n(\gamma/2)}{\Gamma(1 + n\gamma/2)}\left[\left(\Delta_{\varepsilon, b}	\right)^{\gamma-\delta} + \left(\Delta_{\varepsilon, \sigma}	\right)^{\gamma-\delta}\right] + $$ $$+C\int_t^s \dfrac{1}{(s-u)^{1 - \gamma/2}}||\left(\tilde p \otimes H^{n} - \tilde p_{\varepsilon} \otimes H_{\varepsilon}^{n}\right)(t,u,x,y)||_{L_1^{\mu}(\mathbb{R}^d, L_1(\mathbb{R}^d))}du.$$
\bigskip

Continuing the descent until the difference of the first order convolutions, we get
\begin{gather}\label{final}
 ||\left(\tilde p \otimes H^{n+1} - \tilde p_{\varepsilon} \otimes H_{\varepsilon}^{n+1}\right)(t,s,x,y)||_{L_1^{\mu}(\mathbb{R}^d, L_1(\mathbb{R}^d))} \leq \\ C^{n+2}(s-t)^{\delta + (n-1)\gamma/2}\dfrac{\Gamma^n(\gamma/2)}{\Gamma(1 + n\gamma/2)}\left[\left(\Delta_{\varepsilon, b}	\right)^{\gamma-\delta} + \left(\Delta_{\varepsilon, \sigma}	\right)^{\gamma-\delta}\right] + \\+ C^{n} \int_{t}^{s} \frac{1}{\left(s-s_{1}\right)^{1-\gamma/2}} \int_{t}^{s_{1}} \frac{1}{\left(s_{1}-s_{2}\right)^{1-\gamma/2}} \ldots \int_{t}^{s_{n-1}} \frac{||\left(\tilde p \otimes H - \tilde p_{\varepsilon} \otimes H_{\varepsilon}\right)(t,s_n,x,z)||_{L_1^{\mu}(\mathbb{R}^d, L_1(\mathbb{R}^d))}}{\left(s_{n-1}-s_{n}\right)^{1-\gamma/2}}d s_{n} \ldots d s_{1} + \\ +\sum_{i=1}^{n-1} \frac{1}{\Gamma\left(1+\frac{(n-i) \gamma}{2}\right)} \int_{t}^{s} \frac{1}{\left(s-s_{1}\right)^{1-\gamma/2}} \int_{t}^{s_{1}} \frac{1}{\left(s_{1}-s_{2}\right)^{1-\gamma/2}} \ldots \\
 \ldots \int_{t}^{s_{i-2}} \frac{1}{\left(s_{i-2}-s_{i-1}\right)^{1-\gamma/2}} \int_{t}^{s_{i-1}} \frac{\left(s_{i}-t\right)^{{\delta}+\frac{(n-i) \gamma}{2}}}{\left(s_{i-1}-s_{i}\right)^{1-\gamma/2}}\left[\left(\Delta_{\varepsilon, b}	\right)^{\gamma-\delta} + \left(\Delta_{\varepsilon, \sigma}	\right)^{\gamma-\delta}\right](t, s_{i}, x) d s_{i} \ldots ds_{1}. 
\end{gather}

\bigskip

Let us consider the perturbation argument $\Delta_{\varepsilon}(t,s)$. On the diagonal we have $$\Delta_{\varepsilon, b}(t,t,x) + \Delta_{\varepsilon, \sigma}(t,t,x) = \left|b - b_{\varepsilon}\right|(t,x) + \left|\sigma - \sigma_{\varepsilon}\right|(t,x).$$

\bigskip

Clearly, \begin{equation}\label{cont}
	\Delta_{\varepsilon}(t,s) \leq M_{\varepsilon, \alpha, \delta} +M^{\mathcal{C}}_{\varepsilon, \alpha, \delta} \leq  \alpha^{\delta-\gamma/2}\bar{M}_{\varepsilon, \alpha, \delta} + T^{\delta - \gamma/2}\bar{M}^{\mathcal{C}}_{\varepsilon, \alpha, \delta}.
\end{equation}

\bigskip

So, substituting the right-hand side of the estimate \eqref{cont}
 instead of the perturbation argument in \eqref{final} and integrating, we get
$$||\left(\tilde p \otimes H^{n+1} - \tilde p_{\varepsilon} \otimes H_{\varepsilon}^{n+1}\right)(t,s,x,y)||_{L_1^{\mu}\left(\mathbb{R}^d, L_1(\mathbb{R}^d)\right)} \leq $$ $$\leq C^{n+2}\left(\alpha^{\delta-\gamma/2}\bar{M}_{\varepsilon, \alpha, \delta} + T^{\delta - \gamma/2}\bar{M}^{\mathcal{C}}_{\varepsilon, \alpha, \delta}\right)(s-t)^{\delta + (n-1)\gamma/2}\dfrac{1}{\Gamma(1 + n\gamma/2)}\sum\limits_{i=1}^{n} \Gamma^{i}(\gamma/2).$$
\bigskip

The claim of Theorem \ref{Main_thm} follows from the asymptotics of the Gamma function.
 \end{proof}
 \subsection{Derivation of the main result in $L_{\infty}$ case}
 
 It turns out that in the case of uniform perturbations, under assumptions \textbf{(A1)}, \textbf{(A2)},  the difference of transition densities $p(t,s,x,y)$ and $p_{\varepsilon}(t,s,x,y)$ admits an estimate in terms of  $L_q^{\mu} - L_p^{\lambda}$ norm. Recall that for arbitrary function $f(x,y)$, where $x,y \in \mathbb{R}^d$, its $L_q^{\mu} - L_p^{\lambda}$ norm with Borel measures $\mu, \lambda$ and $p, q \in [1, \infty]$ has the following form:
 $$||f(x, y)||_{L_q^{\mu}\left(\lambda(\mathbb{R}^d), L_p\right)} = \left[\int_{\mathbb{R}^d}\left(\int_{\mathbb{R}^d}|f(x, y)|^p \lambda(dy)\right)^{q/p}\mu(dx)\right]^{1/q}.  $$
 \bigskip

For a given parameter $\varepsilon > 0$, let us introduce the perturbation arguments on the fixed time interval $0 \leq t < s \leq T$:
\begin{equation}
	\Delta_{\varepsilon, b}^{\infty} = \sup\limits_{(u, x) \in [t, s] \times \mathbb{R}^d} |b(u,x) - b_{\varepsilon}(u,x)|,
\end{equation}
\begin{equation}
	\Delta_{\varepsilon, \sigma}^{\infty}:=\sup _{u \in [t, s]}\left|\sigma(u, .)-\sigma_{\varepsilon}(u, .)\right|_{\gamma},
	\label{unif_pert}
\end{equation}
where $|\cdot|_{\alpha}$ stands for the usual H\"older norm with exponent $\alpha$ in the space of H\"older continuous bounded functions, i.e.,
$$|f|_{\alpha}:=\sup _{x \in \mathbb{R}^d}|f(t, x)|+\sup _{x \neq y, \ x, y \in \mathbb{R}^{d}} \frac{|f(t,x)-f(t, y)|}{|x-y|^{\alpha}}.$$

 \begin{theorem}
 Let $p(t,s,x,y)$ and $p_{\varepsilon}(t,s,x,y)$ be the transition densities of diffusions \eqref{SDE_0} and \eqref{SDE_1} under \textbf{(A1)}, \textbf{(A2)}.  Then there exists a constant $C > 0$ depending only on $T$ and parameters from the assumptions \textbf{(A)} such that $\forall (t, x), (s, y) \in [0, T] \times \mathbb{R}^d$,
\begin{equation}\label{Main_Est_unif}
	\left|p-p_{\varepsilon}\right|(t, s, x, y) \leq C \cdot \left(\Delta_{\varepsilon}^{\infty}\right)^{\gamma} \cdot \bar{p}(t, s, x, y).
\end{equation}

As a consequence, 
\begin{equation}
||\left(p-p_{\varepsilon}\right)(t, s, x, y)||_{L_q^{\mu}(\lambda({\mathbb{R}^d}), L_p)} \leq C \left(\Delta_{\varepsilon}^{\infty}\right)^{\gamma} \cdot ||\bar{p}(t, s, x, y)||_{L_q^{\mu}(\lambda({\mathbb{R}^d}), L_p)}.
\end{equation}

\end{theorem}
\bigskip
\begin{remark}
	Despite the fact that we defined the uniform perturbations \eqref{unif_pert} of the coefficients taking supremum in spatial variable over the whole $\mathbb{R}^d$ can be reformulated through the flow \eqref{flow}. Indeed, using the semigroup property \eqref{Semigroup_property} of the flow $\theta_{t,s}(\cdot)$, we get
	 $$\sup\limits_{(u, y) \in [t, s] \times \mathbb{R}^d} \left|\left(b-b_{\varepsilon}\right)\left(u, \theta_{u, s}(y)\right)\right| \leq \sup\limits_{(u, x) \in [t, s] \times \mathbb{R}^d} |b(u,x) - b_{\varepsilon}(u,x)| =$$ $$= \sup\limits_{(u, x) \in [t, s] \times \mathbb{R}^d} |b(u,\theta_{u,s}(\theta_{s,u}(x))) - b_{\varepsilon}(u,\theta_{u,s}(\theta_{s,u}(x)))| \leq \sup\limits_{(u, y) \in [t, s] \times \mathbb{R}^d} \left|\left(b-b_{\varepsilon}\right)\left(u, \theta_{u, s}(y)\right)\right|.$$
\end{remark}
\bigskip
\begin{proposition}\label{diff_flows_unif}
	 For all $0 \leq t < s \leq T$ and $y \in \mathbb{R}^d$, there exists a constant $C > 0$ such that
	\begin{equation}
		\left|\theta_{t, s}(y) - \theta^{\varepsilon}_{t, s}(y)\right| \leq C \Delta^{\infty}_{\varepsilon, b} \cdot(s-t). 
	\end{equation}
\end{proposition}
\proof
\begin{gather}
\left|\theta_{t, s}(y) - \theta^{\varepsilon}_{t, s}(y)\right| = \left|\int\limits_{t}^{s}  b(u, \theta_{u, s}(y)) - b_{\varepsilon}(u, \theta^{\varepsilon}_{u, s}(y)) \ du \right| = \\  = \int\limits_{t}^{s} \left| b(u, \theta_{u, s}(y)) - b_{\varepsilon}(u, \theta_{u, s}(y)) +  b_{\varepsilon}(u, \theta_{u, s}(y)) -  b_{\varepsilon}(u, \theta^{\varepsilon}_{u, s}(y)) \right| du \leq \\ \leq \int\limits_{t}^{s} \left| b(u, \theta_{u, s}(y)) - b_{\varepsilon}(u, \theta_{u, s}(y)) \ du \right| + \int\limits_{t}^{s} \left| b_{\varepsilon}(u, \theta_{u, s}(y)) -  b_{\varepsilon}(u, \theta^{\varepsilon}_{u, s}(y)) \ du \right| \leq \\ \leq \Delta_{\varepsilon, b}^{\infty}(s - t) + K \int\limits_{t}^{s} \left|\theta_{u, s}(y) - \theta^{\varepsilon}_{u, s}(y)  \right| \ du.
\end{gather}
By the Gr\"onwall inequality, 
\begin{equation}
\left|\theta_{t, s}(y) - \theta^{\varepsilon}_{t, s}(y)\right| \leq C \cdot \Delta_{\varepsilon, b}^{\infty} \cdot(s-t) .
\end{equation}
\endproof
\bigskip

The next two statements can be easily deduced from the proofs of Lemmas \ref{diff_main_terms} and \ref{diff_par_ker}, respectively.

\bigskip

\begin{lemma}\label{diff_main_terms_unif}
	There exists $C>0$ such that for all $0 \leq t < s \leq T$ , $x, y \in \mathbb{R}^d$, and multiindex $\nu, \ |\nu| \leq 2$:
	\begin{equation}
		\left|D_x^{\nu}\tilde{p}(t,s,x,y) - D_x^{\nu}\tilde{p}_{\varepsilon}(t,s,x,y)\right| \leq \dfrac{C\bar{p}(t,s,x,y)}{(s - t)^{|\nu|/2}} \left(\Delta_{\varepsilon}^{\infty}\right)^{\gamma}.
	\end{equation}

	Consequently, in terms of $L_q^{\mu} - L_p^{\lambda}$ norm, the difference of the main terms can be rewritten as 
	$$\begin{aligned} \label{L_p_diff_main_terms_unif}
		&||\tilde{p}_{\varepsilon} - \tilde{p}||_{L^q_{\mu}(\lambda({\mathbb{R}^d}), L^p)}(t,s,x,y) \leq \\ &\leq \dfrac{C}{(s - t)^{|\nu|/2}}||\bar{p}(t,s,x,y)||_{L^q_{\mu}(\lambda({\mathbb{R}^d}), L^{p})}\left(\Delta_{\varepsilon}^{\infty}\right)^{\gamma}.\end{aligned}$$
\bigskip
\end{lemma}
\bigskip
\begin{lemma}\label{diff_par_ker_unif}
	For all $0 \leq t < s \leq T$ and $x, y \in \mathbb{R}^d$, there exists $C>0$ such that
	$$|H - H_{\varepsilon}|(t,s,x,y) \leq \dfrac{C\bar{p}(t,s,x,y)}{(s-t)^{1 - \gamma/2}} \left(\Delta_{\varepsilon}^{\infty}\right)^{\gamma}.$$
\end{lemma}
\bigskip

To complete the proof of Theorem \ref{Main_Est_unif}, we consider the difference between $n$-th order convolutions. We will derive the respective upper bound by induction.

\begin{lemma}\label{diff_n_order_unif}
	For all $0 \leq t < s \leq T$ and $x, y \in \mathbb{R}^d$, there exists $C>0$ such that
	$$\left|\tilde p \otimes H^{n} - \tilde p_{\varepsilon} \otimes H_{\varepsilon}^{n}\right|(t,s,x,y) \leq (n+1)C^{n + 2} \cdot  \left(\Delta_{\varepsilon}^{\infty}\right)^{\gamma} \cdot (s-t)^{n\gamma/2} \frac{\Gamma^n\left(\gamma/2\right)}{\Gamma(1+n\gamma/2)}\bar{p}(t, s, x, y).$$
	
\end{lemma}
\proof
Let us write decomposition \eqref{decomp}:
\begin{gather}
	\left|\tilde p \otimes H^{n+1} - \tilde p_{\varepsilon} \otimes H_{\varepsilon}^{n+1}\right|(t,s,x,y) \leq \\ \leq \left|\left(\tilde p \otimes H^{n} - \tilde p_{\varepsilon} \otimes H_{\varepsilon}^{n}\right) \otimes H\right|(t,s,x,y) + \left|\left(\tilde p_{\varepsilon} \otimes H_{\varepsilon}^{n}\right) \otimes \left(H - H_{\varepsilon}\right)\right|(t,s,x,y) := \romannumeral 1) + \romannumeral 2). 
\end{gather}
From \eqref{conv_n-order_est} and Lemma \ref{diff_par_ker_unif},
$$\romannumeral 2) \leq  	 C^{n + 2}\left(\Delta_{\varepsilon}^{\infty}\right)^{\gamma} \dfrac{\Gamma\left(\frac{\gamma}{2}\right)^{n}}{\Gamma(\frac{n\gamma}{2})} \int_t^s du \int_{\mathbb{R}^d}dz \dfrac{(u-t)^{n\gamma/2}}{(s - u)^{1 - \gamma/2}} \bar{p}(t, u, x, z)\bar{p}(u,s,z,y) = $$ $$=  C^{n + 2}\left(\Delta_{\varepsilon}^{\infty}\right)^{\gamma} \dfrac{\Gamma\left(\frac{\gamma}{2}\right)^{n + 1}}{\Gamma(1 + \frac{(n+1)\gamma}{2})}(s-t)^{(n+1)\gamma/2}\bar{p}(t,s,x,y). $$ 

Now let $n = 0$. Then 

$$\romannumeral 1) = \left|\left(\tilde p - \tilde p_{\varepsilon}\right) \otimes H\right|(t,s,x,y) \leq C^2 \bar{p}(t,s,x,y) \left(\Delta_{\varepsilon}^{\infty}\right)^{\gamma} \int_t^s (s-u)^{\gamma/2 - 1} du$$ $$ = C^2 \bar{p}(t,s,x,y) \left(\Delta_{\varepsilon}^{\infty}\right)^{\gamma}\dfrac{\Gamma\left(\frac{\gamma}{2}\right)}{\Gamma(1 + \frac{\gamma}{2})}(s-t)^{\gamma/2}.$$

Hence, $$\left|\tilde p \otimes H - \tilde p_{\varepsilon} \otimes H_{\varepsilon}\right|(t,s,x,y) \leq 2C^2 \bar{p}(t,s,x,y) \left(\Delta_{\varepsilon}^{\infty}\right)^{\gamma}\dfrac{\Gamma\left(\frac{\gamma}{2}\right)}{\Gamma(1 + \frac{\gamma}{2})}(s-t)^{\gamma/2}.$$

The induction on $n$ completes the proof.
\endproof

Theorem \ref{Main_Est_unif} now readily follows from Lemmas \ref{diff_main_terms_unif}, \ref{diff_n_order_unif} combined with parametrix expansion \eqref{diff_series} and Stirling's formula.
\subsection{Simulation study}
The aim of this section is to give an example of perturbations with $\Delta_{\varepsilon}(t,s) \rightarrow 0, \ t<s$ when $\varepsilon \rightarrow 0$, but $\left|b(t,x) - b_{\varepsilon}(t,x)\right| \not \rightarrow 0$ uniformly. The idea is to consider the case of highly oscillating
perturbations when the integral is small for  $\varepsilon \rightarrow 0$ , but the difference $\left|b(t,x) - b_{\varepsilon}(t,x)\right|$ remains of a constant order in the
neighborhood of the point $t=\sqrt{\varepsilon}$. We consider the simplest case $d = \gamma = 1$, $\Delta_{\varepsilon, \sigma}(t,s) \equiv 0$, and the drift perturbation of the following form:
$$b(t,x) = x,$$
$$
b_{\varepsilon}(t, x)=b(t, x)+e^{-\frac{t^{2}}{q \varepsilon}} \sin ^{\frac{2}{q}}\left(\frac{t}{\sqrt{\varepsilon}}\right) \cos(x).
$$

Here, $q = \dfrac{p}{p-1}, \ 2 < q < \infty$. Clearly, in this case, \textbf{(A)} holds true. Let us estimate $\Delta_{\varepsilon, b}(t,s)$ when $t < s$.

$$\displaystyle \left|b(t, x)-b_{\varepsilon}(t, x)\right|_{1} \leq  e^{-\frac{t^{2}}{q \varepsilon}} \sin ^{\frac{2}{q}}\left(\frac{t}{\sqrt{\varepsilon}}\right), $$
$$\Delta_{\varepsilon, b}(t, s) \leq \int_{t}^{s} \mathfrak{B}\left(u ; 1, \frac{1}{2}\right) \int_{\mathbb{R}^{d}} \int_{\mathbb{R}^{d}}\left|b\left(u, \theta_{u, s}(y)\right)-b_{\varepsilon}\left(u, \theta_{u, s}(y)\right)\right|_{1} \times $$ $$\times
\bar{p}(t, s, x, y) d y d \mu(x) d u \leq M \int_{t}^{s} \mathfrak{B}\left(u ; 1, \frac{1}{2}\right) e^{-\frac{u^{2}}{q \varepsilon}} \sin ^{\frac{2}{q}}\left(\frac{u}{\sqrt{\varepsilon}}\right) d u \leq $$ $$ 
\leq M\left(\int_{t}^{s} \mathfrak{B}^{p}\left(u ; 1, \frac{1}{2}\right) d u\right)^{\frac{1}{p}}\left(\int_{t}^{s} e^{-\frac{u^{2}}{\varepsilon}} \sin ^{2}\left(\frac{u}{\sqrt{\varepsilon}}\right) d u\right)^{\frac{1}{q}}, \frac{1}{p}+\frac{1}{q}=1. 
$$

First, $$
\left(\int_{t}^{s} \mathfrak{B}^{p}\left(u ; 1, \frac{1}{2}\right) d u\right)^{1 / p}=\frac{1}{B\left(1, \frac{1}{2}\right)(s-t)^{\frac{1}{2}}}\left(\int_{t}^{s}(s-u)^{-\frac{p}{2}} d u\right)^{1 / p} =$$ $$
=\frac{(s-t)^{\frac{1}{p}-1}}{B\left(1, \frac{1}{2}\right)\left(1-\frac{p}{2}\right)^{\frac{1}{p}}}.
$$
\bigskip

Next, using the representation $\sin(x) = \dfrac{\exp(ix) - \exp(-ix)}{2i}$, we get:

$$\int_{0}^{\infty}e^{-x^2}\sin^2(x)dx = \int_{0}^{\infty} \left(\dfrac{e^{-x^2}}{2} - \dfrac{1}{4}e^{-(x-i)^2-1} - \dfrac{1}{4}e^{-(x+i)^2-1}\right)dx = \dfrac{\sqrt{\pi}}{4}\left(1 - 1/e\right).$$

As a consequence, 
$$
\left(\int_{t}^{s} e^{-\frac{u^{2}}{\varepsilon}} \sin ^{2}\left(\frac{u}{\sqrt{\varepsilon}}\right) d u\right)^{\frac{1}{q}} \leq\left(\int_{0}^{\infty} e^{-\frac{u^{2}}{\varepsilon}} \sin ^{2}\left(\frac{u}{\sqrt{\varepsilon}}\right) d u\right)^{\frac{1}{q}}= $$ $$=
\left(\frac{1}{4} \sqrt{\pi \varepsilon}\left(1-1/e\right)\right)^{\frac{1}{q}}.
$$

\bigskip

By Theorem \ref{Main_thm},
$$
\left\|\left(p-p_{\varepsilon}\right)(t, s, x, y)\right\|_{L_{1}^{\mu}\left(\mathbb{R}, L_{1}(\mathbb{R})\right)} \leq $$ $$\leq
\frac{C}{\left(1-\frac{p}{2}\right)^{\frac{1}{p}}\left(\delta-\frac{1}{2}\right)}\left(\varepsilon^{k\left(\delta-\frac{1}{2}\right)}+\varepsilon^{\left(\frac{1}{2}-\kappa\right)\left(1-\frac{1}{p}\right)(1-\delta)}\right),$$
where $0<\kappa<\frac{1}{2}, \frac{1}{2}<\delta<1,1<p<2.
$
\bigskip

Equating the powers, we get for $\delta=\frac{3}{4}, q=\frac{p-1}{p}>2$
$$
\left\|\left(p-p_{\varepsilon}\right)(t, s, x, y)\right\|_{L_{1}^{\mu}\left(\mathbb{R}, L_{1}(\mathbb{R})\right)}  \leq $$ $$
\left.\frac{C}{\left(1-\frac{p}{2}\right)^{\frac{1}{p}}\left(\delta-\frac{1}{2}\right)} \varepsilon^{\frac{(1-\delta)\left(\delta-\frac{1}{2}\right)}{2\left[q\left(\delta-\frac{1}{2}\right)+1-\delta\right]}}\right|_{\delta=\frac{3}{4}}=
\frac{C}{\left(1-\frac{p}{2}\right)^{\frac{1}{p}}} \varepsilon^{1/(8(1+q))}=\frac{C}{\left(1-\frac{p}{2}\right)^{\frac{1}{p}}} \varepsilon^{\frac{1}{24}-\rho}
$$
for any $\rho > 0$ and an appropriate choice of $q \rightarrow 2+$.

Importantly, the perturbation does not tend to zero uniformly. Indeed, for $\frac{\pi \sqrt{\varepsilon}}{6} \leq t \leq \frac{\pi \sqrt{\varepsilon}}{2}$,
$$\left|b_{\varepsilon}(t,x) - b(t,x)\right| \geq 2^{-\frac{2}{q}} e^{-\frac{\pi^{2}}{4 q}} \left|\cos(x)\right|.$$

The following numerical experiments demonstrate the behavior of the perturbation argument $\Delta_{\varepsilon, b}(t,s)$ in the neighborhood of the point $t = \sqrt{\varepsilon}$ for different $\varepsilon$. We consider the time step $\Delta t$, the constant $q=2.01$ and the probability measure $\mu(\cdot) = \delta_1(\cdot)$, where $\delta_{x}(\cdot)$ is a Dirac delta measure centred on some fixed point $x \in \mathbb{R}$.
\begin{figure}[h!]
	\begin{minipage}{0.52\textwidth}
		\centering
		\includegraphics[width=.9\linewidth]{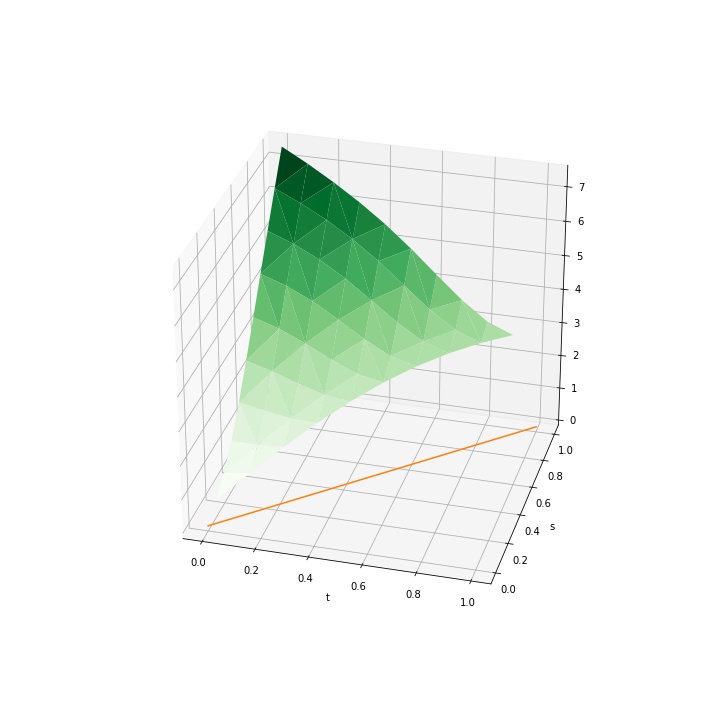}
		\captionsetup{justification=centering}
		\caption{$\Delta_{\varepsilon, b}(t,s), \ 0 \leq t < s \leq 1,$\\\hspace{\textwidth} where  
			$\Delta t = 0.1, \ \varepsilon=1.$}\label{Fig:Data1}
	\end{minipage}\hfill
	\begin{minipage}{0.52\textwidth}
		\centering
		\includegraphics[width=.9\linewidth]{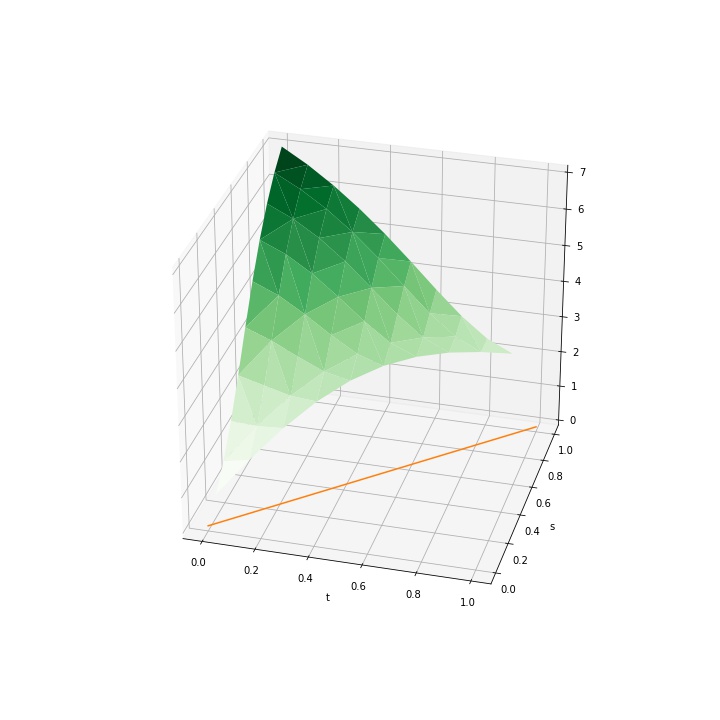}
		\captionsetup{justification=centering}
		\caption{$\Delta_{\varepsilon, b}(t,s), \ 0 \leq t < s \leq 1,$\\\hspace{\textwidth} where $\Delta t = 0.1, \ \varepsilon=0.5.$}\label{Fig:Data2}
	\end{minipage}
\end{figure}
\begin{figure}[h!]
	\begin{minipage}{0.52\textwidth}
		\centering
		\includegraphics[width=.9\linewidth]{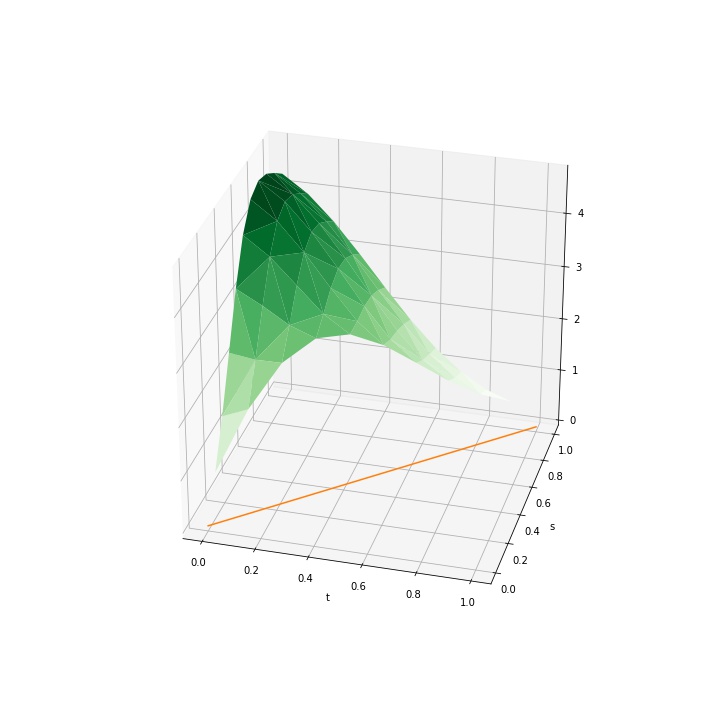}
		\captionsetup{justification=centering}
		\caption{$\Delta_{\varepsilon, b}(t,s), \ 0 \leq t < s \leq 1,$\\\hspace{\textwidth} where $\Delta t = 0.1, \  \varepsilon=0.2.$}\label{Fig:Data3}
	\end{minipage}\hfill
	\begin{minipage}{0.52\textwidth}
		\centering
		\includegraphics[width=.9\linewidth]{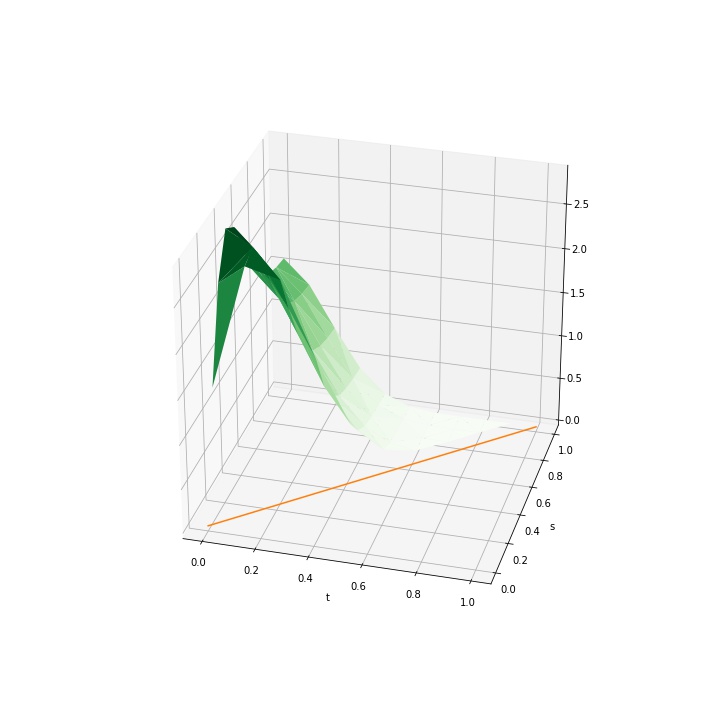}
		\captionsetup{justification=centering}
		\caption{$\Delta_{\varepsilon, b}(t,s), \ 0 \leq t < s \leq 1,$\\\hspace{\textwidth} where $\Delta t = 0.1, \  \varepsilon=0.05.$}\label{Fig:Data4}
	\end{minipage}
\end{figure}

\begin{figure}[h!]
	\begin{minipage}{0.52\textwidth}
		\centering
		\includegraphics[width=.9\linewidth]{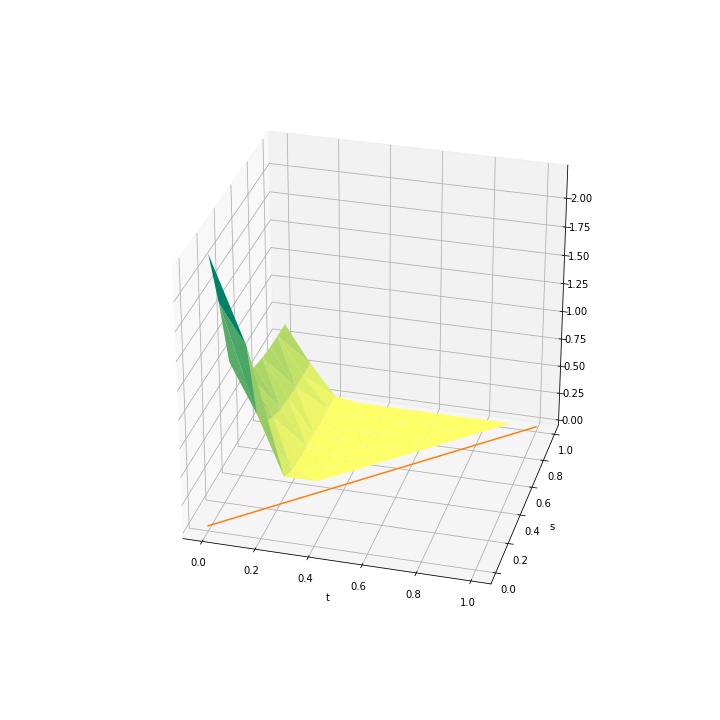}
		\captionsetup{justification=centering}
		\caption{$\Delta_{\varepsilon, b}(t,s), \ 0 \leq t < s \leq 1,$\\\hspace{\textwidth} where $\Delta t = 0.1, \  \varepsilon=0.01.$}\label{Fig:Data5}
	\end{minipage}\hfill
	\begin{minipage}{0.52\textwidth}
		\centering
		\includegraphics[width=.9\linewidth]{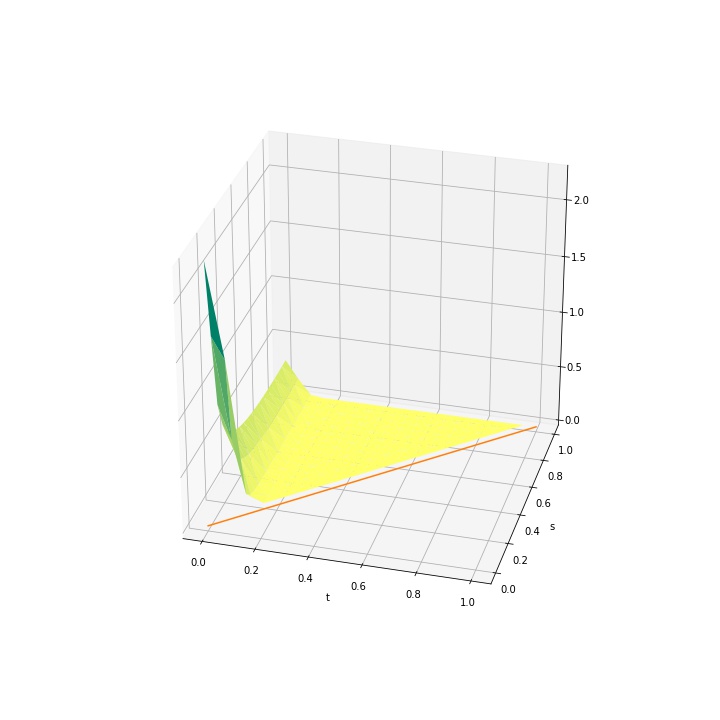}
		\captionsetup{justification=centering}
		\caption{$\Delta_{\varepsilon, b}(t,s), \ 0 \leq t < s \leq 1,$\\\hspace{\textwidth} where $\Delta t = 0.05, \  \varepsilon=0.0025.$}\label{Fig:Data6}
	\end{minipage}
\end{figure}
\pagebreak
\section{Appendix}
\subsection{Proof of Lemma  3.1}
Let us consider the case when $|\nu| = 0$. Let us now identify the transition densities $\tilde{p}(t,s,x,y)$ and $\tilde{p}_{\varepsilon}(t,s,x,y)$ with $(d + 1) \times d$ matrices $\Omega$ and $\Omega^{\varepsilon}$ consisting of the first rows that are the components of the respective mean vectors $\theta_{t,s}(y) = \left(\theta^1_{t,s}(y), \dots, \theta^d_{t,s}(y)\right)$ and $\theta^{\varepsilon}_{t,s}(y)= \left(\theta^{1, \varepsilon}_{t,s}(y), \dots, \theta^{d, \varepsilon}_{t,s}(y)\right)$, and the $d \times d$ covariance matrices, namely,  $\Omega_{ij} = \int_t^s a_{ij}(u, \theta_{u,s}(y))du, \ \Omega^{\varepsilon}_{ij} = \int_t^s a_{\varepsilon}^{ij}(u,\theta^{\varepsilon}_{u,s}(y) )du$, where $i > 1$.
	
\bigskip
For $(d + 1) \times d$ matrix $A$, we will denote by $A_1$ the first row and by $A_{2 : d + 1}$ the square matrix composed of rows from $2$ to $d + 1$.
\bigskip
We can rewrite $\tilde{p}(t,s,x,y)$ and $\tilde{p}_{\varepsilon}(t,s,x,y)$ in terms of $\Omega$ and $\Omega^{\varepsilon}$:

$$\begin{aligned}
	&\tilde{p}(t,s,x,y) = f\left(\Omega\right), \\
 	&\tilde{p}_{\varepsilon}(t,s,x,y) = f\left(\Omega^{\varepsilon}\right),
\end{aligned}$$
	where $$\begin{aligned}
		&f: \mathbb{R}^{(d + 1) \times d} \rightarrow \mathbb{R}, \\
		A \mapsto f(A) = \dfrac{1}{(2\pi)^{d/2}\operatorname{det}(A_{2 : d + 1})^{1/2}} &\exp \left(-\dfrac{1}{2} \langle (A_{2 : d + 1})^{-1}(A_1 - x), A_1 - x \rangle\right).
	\end{aligned}$$

\bigskip

So, applying the Taylor expansion and \eqref{est_der_dens} combined with  \textbf{(A4)} gives:
$$|\tilde{p}_{\varepsilon} - \tilde{p}|(t,s,x,y) = \left|f(\Omega) - f(\Omega_{\varepsilon})\right| = $$ $$ = \left| \sum\limits_{|\nu|=1}\left(\Omega_{\varepsilon}-\Omega\right)^{\nu} \cdot\int\limits_{0}^{1}(1-\lambda)\mathcal{D}^{\nu}f\left\{\Omega + \lambda(\Omega_{\varepsilon}-\Omega)\right\}d\lambda \right| = $$ $$= \left|\sum\limits_{i=1}^d \left(\theta^{i, \varepsilon}_{t,s}(y) - \theta^{i}_{t,s}(y)\right)\cdot\int\limits_{0}^{1}(1-\lambda)\mathcal{D}^{\nu_{1}^i}f\left\{\Omega + \lambda(\Omega_{\varepsilon}-\Omega)\right\}d\lambda  + \right.$$ $$\left.+ \sum\limits_{i,j=1}^d \left(\int_t^s a^{\varepsilon}_{ij}(u,\theta^{\varepsilon}_{u,s}(y) )du - \int_t^s a_{ij}(u,\theta^{\varepsilon}_{u,s}(y) )du\right)\cdot\int\limits_{0}^{1}(1-\lambda)\mathcal{D}^{\nu_{2:d+1}^{ij}}f\left\{\Omega + \lambda(\Omega_{\varepsilon}-\Omega)\right\}d\lambda\right| \leq $$ $$\leq \dfrac{C}{\left({s-t}\right)^{(1 + d)/2}} \exp \left(\dfrac{-\left|\theta_{t,s}(y) + c(\theta_{t,s}^{\varepsilon}(y) - \theta_{t,s}(y)) - x\right|^2}{C(s-t)}\right) \cdot \left(\left|\theta_{t,s}^{\varepsilon}(y) - \theta_{t,s}(y)\right|\right) + $$ $$+ \dfrac{C}{\left({s-t}\right)^{1 + d/2}} \exp \left(\dfrac{-\left|\theta_{t,s}(y) + c(\theta_{t,s}^{\varepsilon}(y) - \theta_{t,s}(y)) - x\right|^2}{C(s-t)}\right) \cdot \left(\int_t^s\left|a(u, \theta_{u,s}(y)) - a_{\varepsilon}(u, \theta^{\varepsilon}_{u,s}(y))\right|du\right) \leq$$  $$\leq \dfrac{C}{\left({s-t}\right)^{1/2}} \bar{p}(t,s,x,y) \exp\left(\frac{c\left|\theta_{t, s}^{\varepsilon}(y)-\theta_{t, s}(y)\right|^{2}}{C(s-t)}\right) \left(\int_{t}^{s}\left|\left(b-b_{\varepsilon}\right)\left(u, \theta_{u, s}(y)\right)\right| d u\right) + $$ $$+\dfrac{C}{\left({s-t}\right)} \bar{p}(t,s,x,y) \exp\left(\frac{c\left|\theta_{t, s}^{\varepsilon}(y)-\theta_{t, s}(y)\right|^{2}}{C(s-t)}\right)\left(
(s-t)\left(\int_{t}^{s}\left|\left(b-b_{\varepsilon}\right)\left(u, \theta_{u, s}(y)\right)\right| d u\right)^{\gamma} +  \right.$$ $$\left. + \int_{t}^{s}\left|\left(\sigma-\sigma_{\varepsilon}\right)\left(u, \theta_{u, s}(y)\right)\right| d u\right) \leq$$ $$\leq  C\bar{p}(t,s,x,y)  \left(
\dfrac{\left(\int_{t}^{s}\left|\left(b-b_{\varepsilon}\right)\left(u, \theta_{u, s}(y)\right)\right| d u\right)^{\gamma}}{(s - t)^{\gamma/2}} +  \dfrac{\left(\int_{t}^{s}\left|\left(\sigma-\sigma_{\varepsilon}\right)\left(u, \theta_{u, s}(y)\right)\right| d u\right)^{\gamma}}{(s-t)^{\gamma}}\right).$$

\bigskip

We will use a slightly different bound:

$$|\tilde{p}_{\varepsilon} - \tilde{p}|(t,s,x,y) \leq $$ \begin{equation}
	 \leq C\bar{p}(t,s,x,y)  \left(
\dfrac{\left(\int_{t}^{s}\left|\left(b-b^{\varepsilon}\right)\left(u, \theta_{u, s}(y)\right)\right| d u\right)^{\gamma - \delta}}{(s - t)^{\gamma/2 - \delta/2}} +  \dfrac{\left(\int_{t}^{s}\left|\left(\sigma-\sigma^{\varepsilon}\right)\left(u, \theta_{u, s}(y)\right)\right| d u\right)^{\gamma - \delta}}{(s-t)^{\gamma - \delta}}\right). \label{est}\end{equation}
\bigskip

Now the upper bound on the difference in $L_1- L_1$ norm easily follows from \eqref{est} and Jensen's inequality:

$$||\left(\tilde{p}_{\varepsilon} - \tilde{p}\right)(t,s,x,y)||_{L_1^{\mu}({\mathbb{R}^d}, L_1(\mathbb{R}^d))} \leq$$ $$\leq \dfrac{C}{(s-t)^{\gamma - \delta}}  \int\limits_{\mathbb{R}^d} \int\limits_{\mathbb{R}^d} \bar{p}(t,s,x,y) \cdot \phi(t,s;y)dy \ \mu(dx) \leq  $$ $$\leq C\left[\left(\int_{\mathbb{R}^d}\int_{\mathbb{R}^d} \int_{t}^{s} \mathfrak{B}(u ; 1, \gamma/2)\left|\left(b-b_{\varepsilon}\right)\left(u, \theta_{u, s}(y)\right)\right|_{1} \bar{p}(t,s,x,y) du  dy \ \mu(dx) \right)^{\gamma - \delta}\right. +$$ $$+\left. \left(\int_{\mathbb{R}^d}\int_{\mathbb{R}^d} \int_{t}^{s} \mathfrak{B}(u ; 1, \gamma/2)\left|\left(\sigma-\sigma_{\varepsilon}\right)\left(u, \theta_{u, s}(y)\right)\right|_{\gamma} \bar{p}(t,s,x,y) du  dy \ \mu(dx)\right)^{\gamma - \delta}\right] = $$ $$=C \left[\left(\Delta_{\varepsilon, b}\right)	^{\gamma-\delta} + \left(\Delta_{\varepsilon, \sigma}\right)	^{\gamma-\delta}\right].$$
\bigskip

The bounds for $|\nu| \geq 1 $ follow from differentiation of the Taylor expansion and similar bounds  \eqref{est_der_dens} for the derivatives of the Gaussian densities $\tilde{p}$ and $\tilde{p}_{\varepsilon}$.
\subsection{Proof of Lemma 3.2}
In order to prove the Lemma \ref{diff_par_ker}, we decompose the difference $\left|H - H_{\varepsilon}\right|(t,s,x,y)$ into 6 parts in the following way:

\begin{gather} 
	\left|H - H_{\varepsilon}\right|(t,s,x,y) \leq \\ \leq \dfrac{1}{2}\sum_{i,j = 1}^d \left|a_{ij}(t,x) - a_{ij}^{\varepsilon}(t,x) - a_{ij}(t,\theta_{t,s}(y)) + a_{ij}^{\varepsilon}(t,\theta_{t,s}(y))\right|\left|\dfrac{\partial^2 \tilde{p}(t,s,x,y)}{\partial x_i \partial x_j}\right| + \\ +\dfrac{1}{2}\sum_{i,j = 1}^d \left|a_{ij}^{\varepsilon}(t,\theta^{\varepsilon}_{t,s}(y)) - a_{ij}^{\varepsilon}(t,\theta_{t,s}(y))\right|\left|\dfrac{\partial^2 \tilde{p}(t,s,x,y)}{\partial x_i \partial x_j}\right| + \\ + \dfrac{1}{2}\sum_{i,j = 1}^d \left|a_{ij}^{\varepsilon}(t,x) - a_{ij}^{\varepsilon}(t,\theta^{\varepsilon}_{t,s}(y))\right|\left|\dfrac{\partial^2 \left(\tilde{p}(t,s,x,y) - \tilde{p}_{\varepsilon}(t,s,x,y)\right)}{\partial x_i  \partial x_j}\right| + \\ + \sum_{i=1}^d \left|b_{i}(t,x) - b_{i}(t, \theta_{t,s}(y))  - b_{i}^{\varepsilon}(t,x) + b_{i}^{\varepsilon}(t, \theta_{t,s}(y))\right| \left|\dfrac{\partial \tilde{p}(t,s,x,y)}{\partial x_i}\right| + \\ + \sum_{i=1}^d \left|b_{i}^{\varepsilon}(t, \theta^{\varepsilon}_{t,s}(y))  -  b_{i}^{\varepsilon}(t, \theta_{t,s}(y))\right|\left|\dfrac{\partial \tilde{p}(t,s,x,y)}{\partial x_i}\right| + \\ + \sum_{i=1}^d \left|b_{i}^{\varepsilon}(t, x)  -  b_{i}^{\varepsilon}(t, \theta^{\varepsilon}_{t,s}(y))\right|\left|\dfrac{\partial \left(\tilde{p} - \tilde{p}_{\varepsilon}\right)(t,s,x,y)}{\partial x_i}\right| := \romannumeral 1) + \romannumeral 2) + \romannumeral 3) + \romannumeral 4) + \romannumeral 5) + \romannumeral 6).
\end{gather}
\bigskip
For $\romannumeral 2)$ and $\romannumeral 5)$, by regularity assumptions and \eqref{est_der_dens}, we have
$$\romannumeral 2) \leq \dfrac{C}{s - t}\bar{p}(t,s,x,y) \left|\theta^{\varepsilon}_{t,s}(y) - \theta_{t,s}(y)\right|^{\gamma} \leq \dfrac{C}{(s-t)^{1 -\delta/2}}\bar{p}(t,s,x,y)\left(\int_t^s \left|\left(b-b^{\varepsilon}\right)\left(u, \theta_{u, s}(y)\right)\right|_{1} d u\right)^{\gamma-\delta}.$$
$$\romannumeral 5) \leq \dfrac{C}{(s-t)^{1/2}}\bar{p}(t,s,x,y) \left|\theta^{\varepsilon}_{t,s}(y) - \theta_{t,s}(y)\right| \leq \dfrac{C}{(s-t)^{1/2}}\bar{p}(t,s,x,y)\left(\int_t^s \left|\left(b-b^{\varepsilon}\right)\left(u, \theta_{u, s}(y)\right)\right|_{1} d u\right).$$
\bigskip

To obtain an upper bound for $\romannumeral 1)$, we first estimate the first part of the corresponding expression:
$$\left|a(t,x) - a^{\varepsilon}(t,x) - a(t,\theta_{t,s}(y)) + a^{\varepsilon}(t,\theta_{t,s}(y))\right| =$$ $$= \left|\sigma \sigma^*(t,x) - \sigma^{\varepsilon}(\sigma^{\varepsilon})^*(t,x) - \sigma \sigma^*(t,\theta_{t,s}(y)) + \sigma^{\varepsilon}(\sigma^{\varepsilon})^*(t,\theta_{t,s}(y))\right| \leq $$ $$\leq \left[ \left|\sigma(t,x) - \sigma^{\varepsilon}(t,x) - \left(\sigma (t,\theta_{t,s}(y)) - \sigma^{\varepsilon}(t,\theta_{t,s}(y))\right)\right|\left|\sigma^*(t,x)\right| + \right.$$ $$+ \left|\sigma(t,x) - \sigma^{\varepsilon}(t,x) - \left(\sigma (t,\theta_{t,s}(y)) - \sigma^{\varepsilon}(t,\theta_{t,s}(y))\right)\right|\left|(\sigma^{\varepsilon})^*(t,x)\right| +$$ $$+ \left|\sigma^*(t,x) - (\sigma^{\varepsilon})^*(t,x) - \left(\sigma^* (t,\theta_{t,s}(y)) - (\sigma^{\varepsilon})^*(t,\theta_{t,s}(y))\right)\right|\left|\sigma(t,\theta_{t,s}(y))\right| + $$ $$\left.+\left|\sigma (t,\theta_{t,s}(y)) - \sigma^{\varepsilon}(t,\theta_{t,s}(y))\right| \left|(\sigma^{\varepsilon})^*(t,x) - (\sigma^{\varepsilon})^*(t,\theta_{t,s}(y))\right| + \right. $$ $$\left.   +\left|\sigma^{\varepsilon} (t,x) - \sigma^{\varepsilon}(t,\theta_{t,s}(y))\right| \left|\sigma^*(t,x) - (\sigma^{\varepsilon})^*(t,x)\right|\right] \leq$$ $$\leq C\left(\left|x - \theta_{t,s}(y)\right|^{\gamma} + \left|x - \theta_{t,s}(y)\right|^{2\gamma} \right)\left|\left(\sigma - \sigma^{\varepsilon}\right)(t, \theta_{t,s}(y))\right|_{\gamma}.$$

Combining the estimate above and Lemma \ref{est_der_dens} with the semigroup property of the flows, we have for $\romannumeral 1)$ and $\romannumeral 4)$
$$\romannumeral 1) \leq  \dfrac{C}{s - t}\bar{p}(t,s,x,y)\left|a(t,x) - a^{\varepsilon}(t,x) - a(t,\theta_{t,s}(y)) + a^{\varepsilon}(t,\theta_{t,s}(y))\right| \leq$$ $$\leq \dfrac{C}{s - t}\bar{p}(t,s,x,y)\left(\left|x - \theta_{t,s}(y)\right|^{\gamma} + \left|x - \theta_{t,s}(y)\right|^{2\gamma} \right)\left|\left(\sigma - \sigma^{\varepsilon}\right)(t, \theta_{t,s}(y))\right|_{\gamma} \leq \dfrac{C\bar{p}(t,s,x,y)}{(s - t)^{1 - \gamma/2}}\left|\left(\sigma - \sigma^{\varepsilon}\right)(t, \theta_{t,s}(y))\right|_{\gamma}.$$

$$\romannumeral 4) \leq \dfrac{C\bar{p}(t,s,x,y)}{(s - t)^{1/2}}\left|\left(b - b^{\varepsilon}\right)(t, \theta_{t,s}(y))\right|_{1}.$$

\bigskip
Finally, using the control Lemma \ref{diff_main_terms} for $\romannumeral 3)$ and $\romannumeral 6)$, we get

$$\romannumeral 3) \leq C\bar{p}(t,s,x,y)  \left(
\dfrac{\left(\int_{t}^{s}\left|\left(b-b^{\varepsilon}\right)\left(u, \theta_{u, s}(y)\right)\right|_{1} d u\right)^{\gamma - \delta}}{(s - t)^{1 - \delta/2}} +  \dfrac{\left(\int_{t}^{s}\left|\left(\sigma-\sigma^{\varepsilon}\right)\left(u, \theta_{u, s}(y)\right)\right|_{\gamma} d u\right)^{\gamma - \delta}}{(s-t)^{1 + \gamma/2 - \delta}}\right).$$

$$\romannumeral 6) \leq \left(
\dfrac{\left(\int_{t}^{s}\left|\left(b-b^{\varepsilon}\right)\left(u, \theta_{u, s}(y)\right)\right|_{1} d u\right)^{\gamma - \delta}}{(s - t)^{\gamma/2 - \delta/2}} +  \dfrac{\left(\int_{t}^{s}\left|\left(\sigma-\sigma^{\varepsilon}\right)\left(u, \theta_{u, s}(y)\right)\right|_{\gamma} d u\right)^{\gamma - \delta}}{(s-t)^{\gamma - \delta}}\right).$$
\bigskip

Summing up the bounds derived above, we complete the proof.

\bibliographystyle{acm}

\end{document}